\documentclass[12pt,leqno]{article}

\usepackage{amsmath,amssymb,amsthm}
\usepackage[all]{xy}
\usepackage{lscape}

\makeatletter

\swapnumbers

\newtheorem{theorem}{Theorem}[section]
\newtheorem{proposition}[theorem]{Proposition}
\newtheorem{lemma}[theorem]{Lemma}
\newtheorem{corollary}[theorem]{Corollary}

\theoremstyle{definition}

\newtheorem{example}[theorem]{Example}
\newtheorem{definition}[theorem]{Definition}

\theoremstyle{remark}

\theoremstyle{remark}
\newtheorem{remark}[theorem]{Remark}

\def\({{\rm (}}
\def\){{\rm )}}

\let\Mathrm\operator@font
\let\Cal\mathcal
\newcommand{\fm}{\ensuremath{\mathfrak m}}
\newcommand{\fn}{\ensuremath{\mathfrak n}}
\newcommand{\fp}{\ensuremath{\mathfrak p}}
\newcommand{\fq}{\ensuremath{\mathfrak q}}
\newcommand{\fa}{\ensuremath{\mathfrak a}}
\renewcommand{\O}{\Cal O}
\newcommand{\cO}{\Cal O}
\newcommand{\I}{\Cal I}
\newcommand{\J}{\Cal J}
\newcommand{\K}{\Cal K}
\newcommand{\cP}{\Cal P}
\renewcommand{\P}{\Cal P}
\newcommand{\Q}{\Cal Q}
\newcommand{\N}{\Cal N}
\newcommand{\M}{\Cal M}
\newcommand{\cL}{\Cal L}
\renewcommand{\L}{\Cal L}

\def\mt#1{$#1$}

\def\standop#1{\mathop{\Mathrm #1}\nolimits}
\def\difstop#1#2{\expandafter\def\csname #1\endcsname{\standop{#2}}}
\def\defstop#1{\difstop{#1}{#1}}

\defstop{AB}
\defstop{ann}
\defstop{Ass}

\defstop{codim}
\defstop{Coh}
\defstop{Coker}
\defstop{Cone}

\defstop{depth}

\defstop{EM}
\defstop{End}
\defstop{ev}
\defstop{Ext}

\defstop{Fitt}
\defstop{Flat}
\defstop{Func}

\defstop{Gr}

\difstop{height}{ht}
\defstop{Hom}

\def\id{\mathord{\Mathrm{id}}}

\difstop{Image}{Im}
\defstop{Im}

\defstop{Ker}
\defstop{Kos}

\defstop{Lch}
\defstop{length}
\defstop{Lin}
\defstop{LQ}
\defstop{Lqc}
\defstop{lqc}

\defstop{Mat}
\defstop{Min}
\defstop{Mod}
\defstop{Mor}

\defstop{PA}
\defstop{PM}
\defstop{Proj}

\defstop{Qch}
\defstop{qch}

\defstop{rank}
\def\red{_{\Mathrm{red}}}
\defstop{res}
\defstop{rad}

\def\Sch{\underline{\Mathrm Sch}}

\defstop{Spec}
\defstop{supp}
\defstop{Supp}
\defstop{Sym}

\difstop{tdeg}{trans.deg}

\defstop{Tor}
\defstop{type}
\difstop{trace}{tr}
\difstop{Sim}{SIm}
\def\GSim{\mathop{G\mathrm{SIm}}\nolimits}

\def\grad#1{\sqrt[G]{#1}}

\let\indlim\varinjlim

\def\uHom{\mathop{\mbox{\underline{$\Mathrm Hom$}}}\nolimits}

\def\uGamma{\mathop{\mbox{\underline{$\Gamma$}}}\nolimits}
\def\uH{\mathop{\mbox{\underline{$H$}}}\nolimits}
\def\uann{\mathop{\mbox{\underline{$\Mathrm ann$}}}\nolimits}
\def\uFitt{\mathop{\mbox{\underline{$\Mathrm Fitt$}}}\nolimits}

\def\sdarrow#1{\downarrow\hbox to 0pt{\scriptsize$#1$\hss}}
\def\suarrow#1{\uparrow\hbox to 0pt{\scriptsize$#1$\hss}}
\def\ssearrow#1{\searrow\hbox to 0pt{\scriptsize$#1$\hss}}

\def\ext{{\textstyle\bigwedge}}

\def\section{\@startsection{section}{1}{\z@ }%
{-3.5ex plus -1ex minus -.2ex}{2.3ex plus .2ex}{\bf }}

\long\def\refname{\par\kern -3ex
\begin{center}\rm R\sc{eferences}\end{center}\par\kern 
-2ex}

\def\@seccntformat#1{\csname the#1\endcsname.\quad}

\def\@@@sect#1#2#3#4#5#6[#7]#8{%
   \ifnum #2>\c@secnumdepth 
      \def \@svsec {}\else \refstepcounter {#1}%
      \def\@svsec{}
   \fi 
   \@tempskipa #5\relax 
   \ifdim \@tempskipa >\z@ 
     \begingroup #6\relax \@hangfrom {\hskip #3\relax 
     \@svsec}{\interlinepenalty \@M #8\par }\endgroup 
     \csname #1mark\endcsname {#7}
   \else 
   \def \@svsechd {#6\hskip #3\@svsec #8\csname #1mark\endcsname {#7}}
   \fi \@xsect {#5}}

\def\@@@startsection#1#2#3#4#5#6{%
 \if@noskipsec \leavevmode \fi \par \@tempskipa #4\relax \@afterindenttrue 
 \ifdim \@tempskipa <\z@ \@tempskipa -\@tempskipa \@afterindentfalse 
 \fi \if@nobreak \everypar {}\else \addpenalty {\@secpenalty }\addvspace 
  {\@tempskipa }\fi \@ifstar {\@ssect {#3}{#4}{#5}{#6}}{\@dblarg 
  {\@@@sect {#1}{#2}{#3}{#4}{#5}{#6}}}}

\def\theparagraph{\thesection.\arabic{paragraph}}
\def\aparagraph{\@@@startsection{paragraph}{2}{\z@ }%
              {1.75ex plus .2ex minus .15ex}{-1em}{\bf(\theparagraph) } }
\def\paragraph{\@@@startsection{paragraph}{2}{\z@ }%
              {1.75ex plus .2ex minus .15ex}{-1em}{}{\bf(\theparagraph)} }

\let\c@theorem\c@paragraph

\title{$G$-prime and $G$-primary $G$-ideals on $G$-schemes}

\author{M{\sc itsuyasu} H{\sc ashimoto}
and 
M{\sc itsuhiro} M{\sc iyazaki}}

\date{\normalsize
Graduate School of Mathematics, Nagoya University\\
Chikusa-ku,  Nagoya 464--8602 JAPAN\\
{\small{\tt hasimoto@math.nagoya-u.ac.jp}}\\
~\\
Department of mathematics, Kyoto University of Education\\
1~Fukakusa-Fujinomori-cho, Fushimi-ku Kyoto 612--8522 JAPAN\\
{\small {\tt g53448@kyokyo-u.ac.jp}}
\\
~
\\
Dedicated to Professor Masayoshi Nagata
}

\begin{document}

\maketitle
\renewcommand\thefootnote{\fnsymbol{footnote}}%
\footnote[0]{\textit{Key Words and Phrases}.
$G$-prime $G$-ideal, $G$-primary $G$-ideal, $G$-associated $G$-prime $G$-ideal, 
Matijevic--Roberts type theorem.
2010 \textit{Mathematics Subject Classification}.
Primary 14L30.}%
\begin{abstract}
Let $G$ be a flat finite-type group scheme over a scheme $S$,
and $X$ a noetherian $S$-scheme on which $G$ acts.
We define and study $G$-prime and $G$-primary $G$-ideals on $X$
and study their basic properties.
In particular, we prove the existence of minimal $G$-primary decomposition
and the well-definedness of $G$-associated $G$-prime $G$-ideals.
We also prove a generalization of Matijevic--Roberts type theorem.
In particular, we prove Matijevic--Roberts type theorem on graded rings
for $F$-regular and $F$-rational properties.
\end{abstract}

\section{Introduction}
In this introduction, let $R$ be a noetherian base ring, $G$ a flat 
group scheme of finite type over $R$, and consider a noetherian 
$R$-algebra $A$ with
a $G$-action, for simplicity.

Let $H$ be a finitely generated 
abelian group, $R=\Bbb Z$, and $W=RH$ the group algebra.
Letting each $h\in H$ group-like, $W$ is a finitely generated
flat commutative Hopf algebra over
$R$, and hence $G=\Spec W$ is an affine flat $R$-group scheme of finite type.
It is well-known that a $G$-algebra is nothing but an $H$-graded ring,
and for a $G$-algebra $A$, a $(G,A)$-module is nothing but a graded $A$-module.
This is the most typical and important case, and actually many of 
our ideas and results in this paper for general $G$
already appeared in \cite{GW1}, \cite{GW2}, \cite{Kamoi}, 
\cite{RA}, 
and \cite{PK}
as those for $H$-graded rings.
Note that in \cite{Kamoi}
and \cite{RA}, 
non-finitely generated $H$ is treated (until section~5, our $G$ need not be
of finite type, and so the case that $H$ is not finitely generated is
also treated in our discussion).

Let $\fp$ be a prime ideal of a $\Bbb Z^n$-graded ring $A$.
Then $\fp^*$, the largest homogeneous ideal contained in $\fp$, is again a
prime ideal.
An associated prime ideal of a homogeneous 
ideal of a noetherian $\Bbb Z^n$-graded ring
is again homogeneous.
These well-known facts on graded rings
can be generalized to results on actions of 
smooth groups with connected fibers, see 
Corollary~\ref{regular-and-connected.thm}.

However, this is not true any more for more general group scheme actions.
For example, these results fail to be true for torsion-graded rings.
But even for general $G$-algebra for a group scheme $G$, 
the ideal of the form $\fp^*$ is very special among other $G$-ideals,
where $\fp^*$ is defined to be the largest $G$-ideal contained in $\fp$.
In general, we define that a $G$-ideal $P$ is {\em $G$-prime}
if $P=\fp^*$ for some prime ideal $\fp$ of $A$.
It is not difficult to show that a $G$-ideal $P$ 
of $A$ is $G$-prime if and only if the following holds:
$P\neq A$, and if $I$ and $J$ are
$G$-ideals and $IJ\subset P$, then either $I\subset P$ or $J\subset P$,
see Lemma~\ref{G-prime.thm}.
Thus our definition is a straightforward generalization of Kamoi's
definition of an $H$-prime homogeneous ideal \cite[Definition~1.2]{Kamoi}.

A $G$-primary $G$-ideal is defined similarly, and our definition is
a natural generalization of that of \cite{RA} and \cite{PK}.
The purpose of this paper is to define $G$-prime and $G$-primary $G$-ideal
and study basic properties of these $G$-ideals.
$G$-radicals of $G$-ideals are also defined and studied.
We also define and study basic properties of $G$-radical and $G$-maximal
$G$-ideals.

One of the most important motivation of our study is a generalization of 
so called Matijevic--Roberts  type theorem.
This type of theorem asserts that 
if $A$ is a noetherian ring with a $G$-action, $P$ a prime ideal of $A$,
and if $A_{P^*}$ enjoys the property $\Bbb P$, then $A_{P}$ has the
same property $\Bbb P$, where $\Bbb P$ is either \lq Cohen--Macaulay,'
\lq Gorenstein,' \lq complete intersection,' or \lq regular.'
The theorem was originally conjectured by Nagata \cite{Nagata2} for 
$\Bbb Z$-graded rings and Cohen--Macaulay property.
The theorem for $\Bbb Z^n$-graded rings was proved by
Hochster--Ratliff \cite{HRa},
Matijevic--Roberts \cite{MR},
Aoyama--Goto \cite{AG1},
Matijevic \cite{Matijevic},
Goto--Watanabe \cite{GW2},
Cavaliere--Niesi \cite{CN},
and
Avramov--Achilles \cite{AA}.
The theorem was then generalized to 
the action of affine smooth $G$ with connected fibers 
\cite[Theorem~II.2.4.2]{Hashimoto}.
The affine assumption was recently removed by M. Ohtani and the first author
(unpublished).
Note that $A_{P^*}$ makes sense in these cases
because $P^*$ is a prime ideal.
On the other hand, Kamoi \cite[Theorem~2.13]{Kamoi} 
proved the theorem for graded rings, graded by
general $H$, for Cohen--Macaulay and Gorenstein properties.
As $P^*$ is not a prime ideal any more, he modified the statement of the
theorem.

Although we need to assume that either $G$ is smooth or
$R=k$ is a perfect field for the property \lq regular,' 
we prove this theorem for general $G$, see Corollary~\ref{M-R.thm}.
Assuming that $A$ is locally excellent
and also assuming that $G$ is smooth or $R=k$ is a perfect field, 
we also prove the theorem for $F$-regularity and $F$-rationality.
It seems that this assertion has not been known even 
as a theorem on $\Bbb Z^n$-graded rings.
Note that $P^*$ is not a prime ideal in the general settings, and
as in Kamoi's work, we need to modify the statement.
That is, we replace the condition \lq $A_{P^*}$ enjoys $\Bbb P$' by
\lq for some minimal prime ideal $\fp$ of $P^*$, $A_{\fp}$ enjoys $\Bbb P$.'
We can prove that ``for some minimal prime ideal'' is equivalent to ``for
any minimal prime ideal,'' see Corollary~\ref{independence2.thm}.

We also discuss $G$-primary decomposition.
It is an analogue of primary decomposition.
Our main results are the existence and the uniqueness of the $G$-primary
decomposition.
The case for $H$-graded rings is treated in \cite{RA} and \cite{PK}.
We also prove that a $G$-primary $G$-ideal does not have an embedded prime 
ideal (see Corollary~\ref{G-primary-no-emb.thm}).
There is a deep connection between primary decompositions of a $G$-ideal $I$
and those of \lq $G$-primary components' of $I$, see Theorem~\ref{main.thm}.
For related results on graded rings, see \cite[Proposition~2.2]{Kamoi}
and \cite[Corollary~4.2]{PK}.

When we consider a group action, considering only affine schemes is
sometimes too restrictive, even if a goal is a theorem on rings.
We treat a group scheme action on a noetherian scheme $X$.

We also prove some scheme theoretic properties on $G$-schemes such that
$0$ is $G$-primary.
If $X$ is a noetherian $G$-scheme such that $0$ is $G$-primary, then
the dimension of the fiber of the second projection $p_2:G\times X\rightarrow
X$ is constant (Proposition~\ref{constant.thm}).
If, moreover, $X$ is of finite type over a field or $\Bbb Z$ (more generally, 
$X$ is {\em Ratliff}, see (\ref{Ratliff.par})), then $X$ is
equidimensional (Proposition~\ref{equidimensional.thm}).
We say that a ring $B$ is {\em Ratliff} if $B$ is noetherian, universally 
catenary, Hilbert, and $B/P$ satisfies the first chain condition for
any minimal prime ideal $P$ of $B$.
A field and $\Bbb Z$ are Ratliff.
A finite-type algebra over a Ratliff ring is Ratliff 
(Lemma~\ref{Ratliff-of-finite-type.thm}).
We say that a scheme $Y$ is Ratliff if $Y$ has a finite open covering
consisting of prime spectra of Ratliff rings.

Section~2 is preliminaries.
We review basics on scheme theoretic image, Fitting ideals, 
and $(G,\O_X)$-modules.
In section~3, we study primary decompositions of ideal sheaves over a
noetherian scheme.
In section~4, we define and study $G$-prime and $G$-radical $G$-ideals.
We also study some basic properties of 
$\fn^*$ (the largest quasi-coherent $(G,\O_X)$-module 
contained in $\fn$) used later.
In section~5, assuming that the $G$-scheme $X$ is noetherian, we 
define and study $G$-primary $G$-ideals on $X$ 
and $G$-primary decomposition.
We prove that $G$-associated $G$-prime $G$-ideal is well-defined.
In section~6, assuming that $G$ is of finite type (more generally, 
the second projection $p_2:G\times X\rightarrow X$ is of finite type),
we study further properties of $G$-prime and $G$-primary $G$-ideals.
In section~7, we
prove a generalization of Matijevic--Roberts type theorem
as explained above.

After the former version of this paper was posted on the arXiv, 
we got aware of the important paper
of
Perling and Trautmann \cite{PT}.
Although they treat only linear algebraic (smooth) groups over
an algebraically closed field, 
there are some overlaps with \cite{PT} and ours.
In particular, \cite[Theorem~4.11]{PT} for $X$ of finite type follows from
our Lemma~\ref{existence.thm}, Corollary~\ref{primary3.thm}, and
Corollary~\ref{regular-and-connected.thm}.
Note that disconnected $G$ is also treated in \cite{PT}.
Our paper is basically written independently of \cite{PT}, but 
we added a paragraph (\ref{PT.par}) which shows
how \cite[Theorem~4.18]{PT} (for the case that $X$ is quasi-compact)
is proved from our approach.
The largest difference of our paper from \cite{PT} is that we 
define and studied $G$-primary $(G,\O_X)$-module for a general group 
scheme.

\section{Preliminaries}

\paragraph
A scheme is not required to be separated in general.

Let $Y$ be a scheme and $\Cal I$ a quasi-coherent ideal of $\Cal O_Y$.
Then we denote the closed subscheme defined by $\Cal I$ by $V(\Cal I)$.
For an ideal $\Cal J$ of $\Cal O_Y$, the sum of all quasi-coherent
ideals of $\Cal O_Y$ contained in $\Cal J$ is denoted by $\Cal J^\star$.
Note that $\Cal J^\star$ is the largest quasi-coherent ideal contained in
$\Cal J$.

Let $\varphi:Y\rightarrow Z$ be a morphism of schemes.
Then $V(\Ker(\Cal O_Z\rightarrow \varphi_*\Cal O_Y)^\star)$ is denoted by 
$\Sim\varphi$, and called the {\em scheme theoretic image} (or the closed 
image) of $\varphi$ \cite[(9.5.3)]{EGA-I}.
Let $\psi:Y\rightarrow \Sim\varphi$ be the induced map.
If $W$ is a separated scheme, then 
\[
\psi^*: \Hom_{\Sch}(\Sim\varphi,W)\rightarrow
\Hom_{\Sch}(Y,W)
\]
is injective, where $\Sch$ denotes the category of schemes 
\cite[(9.5.6)]{EGA-I}.

\paragraph\label{reduced-sim.par}
By definition, $\Sim\varphi$ is the smallest closed subscheme
of $Z$ through which $\varphi$ factors.
So 
it is easy to see that
if $Y$ is reduced, then the closure of 
the image $\varphi(Y)$ of $Y$ by $\varphi$ with the reduced structure is
$\Sim\varphi$ \cite[Exercise~II.3.11 (d)]{Hartshorne}.
In particular, $\Sim\varphi$ is reduced if $Y$ is reduced.

\begin{example}
Let $f:A\rightarrow B$ be a ring homomorphism.
Then the scheme theoretic image of ${}^af:\Spec B\rightarrow \Spec A$ is
$\Spec (A/\Ker f)\cong\Spec \Image f$.
\end{example}

\paragraph Let $Y$ be a scheme, and $Z$ a subscheme of $Y$.
The scheme theoretic image $\Sim\iota$ is called the {\em closure} 
of $Z$ in $Y$, and is denoted by $\bar Z$, where $\iota:Z\hookrightarrow Y$ 
is the inclusion.

\paragraph Let $\varphi:Y\rightarrow Z$ be a faithfully flat morphism.
Then for an $\Cal O_Z$-module $\Cal M$, the unit of adjunction
$u:\varphi^*\Cal M\rightarrow \varphi^*\varphi_*\varphi^*\Cal M$ is
a split mono.
Since $\varphi^*$ is faithful exact, $u:\Cal M\rightarrow \varphi_*\varphi^*
\Cal M$ is also monic.
In particular, $\Cal O_Z\rightarrow \varphi_*\Cal O_Y$ is a mono, and
hence $\Sim\varphi=Z$.

\paragraph\label{sequence-sim.thm}
Let
\[
Y\xrightarrow\varphi Z \xrightarrow\psi W
\]
be a sequence of morphisms of schemes.
Let $\iota:\Sim\varphi\rightarrow Z$ be the inclusion.
Then $\Sim(\psi\varphi)=\Sim(\psi\iota)$ \cite[(9.5.5)]{EGA-I}.

In particular, if $\varphi:Y\rightarrow Z$ is a flat finite-type morphism
between noetherian schemes, then $\Sim\varphi=\overline{\Im\varphi}$.

\paragraph 
If $\varphi:Y\rightarrow Z$ is a quasi-compact morphism, then
$\Ker(\Cal O_Z\rightarrow\varphi_*\Cal O_Y)$ is quasi-coherent, and hence
$\Sim(\varphi)=V(\Ker(\Cal O_Z\rightarrow\varphi_*\Cal O_Y))$.
If this is the case, 
$\Sim(\varphi)$ agrees with the closure of the image $\varphi(Y)$
of $\varphi$, set-theoretically.

To verify this, we may assume that $Z=\Spec A$ is affine, and hence $Y$ is
quasi-compact.
There is a finite affine open covering $(U_i)$ of $Y$.
Let $U:=\coprod_i U_i \rightarrow Y$ be the canonical map.
Replacing $Y$ by $U$, we may assume that $Y=\Spec B$ is also affine.
Then $\varphi_*\Cal O_Y$ is quasi-coherent, and hence
$\Ker(\Cal O_Z\rightarrow\varphi_*\Cal O_Y)$ is also quasi-coherent.
It remains to show that if $A\hookrightarrow B$ is an inclusion of a
subring, then the associated map $\varphi:
\Spec B\rightarrow \Spec A$ is dominating.
Note that for $P\in\Spec A$, there exists some minimal prime ideal 
$Q\in\Spec A$
such that $Q\subset P$ (this follows from Zorn's lemma).
As $0\neq A_Q\rightarrow B_Q=\kappa(Q)\otimes_A B
$ is injective, $\Spec (\kappa(Q)\otimes_A B)$ is non-empty.
This shows $Q\in \varphi(\Spec B)$, and we are done.

Note that in general, if $\varphi$ is quasi-compact quasi-separated, 
then $\varphi_*\Cal O_Y$ is quasi-coherent \cite[(9.2.1)]{EGA-I}.

If $\varphi$ is not quasi-compact, then $\Sim(\varphi)$ need not agree with
the closure of $\varphi(Y)$.
For example, let $Z=\Spec R$, where $R$ is a DVR with a uniformizing 
parameter $t$.
Let $Y=\coprod_{i\geq 1}\Spec R/(t^i)$, and $\varphi:Y\rightarrow Z$ be
the canonical map.
Then the closure of the image of $\varphi$ is the closed point of $Z$
set-theoretically, but $\Sim(\varphi)$ is the whole $Z$.

\paragraph \label{sim-flat.par}
(cf.\ \cite[(9.5.8)]{EGA-I})
Let 
\[
\xymatrix{
Y' \ar[r]^{\varphi'} \ar[d]^g & Z' \ar[d]^f \\
Y \ar[r]^\varphi & Z
}
\]
be a cartesian square of schemes.
Assume that $\varphi$ is quasi-compact quasi-separated, and 
$f$ is flat.
Then $\Sim\varphi'=f^{-1}(\Sim\varphi)$.
Indeed, if 
\[
0\rightarrow \Cal I \rightarrow \Cal O_Z \rightarrow \varphi_*\Cal O_Y
\]
is exact, then applying the exact functor $f^*$ to it,
\[
0\rightarrow f^*\Cal I\rightarrow \Cal O_{Z'} \rightarrow
(\varphi')_*\Cal O_{Y'}
\]
is exact.

\paragraph
Let $Y$ be a scheme.
For a closed subscheme $W$ of $Y$, we denote the defining ideal sheaf 
of $W$ by $\Cal I(W)$.
For a 
morphism $f:Y\rightarrow Z$ of schemes and an ideal sheaf
$\Cal I$ of $\Cal O_Y$, the ideal sheaf $(\eta^{-1}(f_*\Cal I))^\star$ is 
denoted by $\Cal I\cap \Cal O_Z$, where $\eta:\Cal O_Z\rightarrow
f_*\Cal O_Y$ is the canonical map.
In other words, 
$\Sim(V(\Cal I)\hookrightarrow Y\rightarrow Z)=V(\Cal I \cap \Cal O_Z)$.

If $Y\xrightarrow{\varphi}Z\xrightarrow{\psi}W$ is a sequence of morphisms
and $\I$ is an ideal of $\Cal O_Y$, then
$(\I\cap \cO_Z)\cap \cO_W=\I\cap O_W$ by 
(\ref{sequence-sim.thm}).

For a morphism $f:Y\rightarrow Z$ and an ideal $\J$ of $\Cal O_Z$, 
the image of $f^*\J\rightarrow f^*\Cal O_Z\rightarrow \Cal O_Y$ is
denoted $\J\Cal O_Y$.
If $\J$ is a quasi-coherent ideal, then so is $\J\Cal O_Y$, and
$V(\J\Cal O_Y)=f^{-1}(V(\J))$.
Note that $\J\subset \J\Cal O_Y\cap\Cal O_Z$ and
$\I\supset(\I\cap\Cal O_Z)\Cal O_Y$.

If $f:Y\rightarrow Z$ is a quasi-compact 
immersion, then $(\Cal I\cap \Cal O_Z)\Cal O_Y$ is
$\Cal I$.

\paragraph Let $R$ be a commutative ring.
Let $\varphi:F\rightarrow E$ be a map of $R$-free modules with $E$ finite,
and $n$ be an integer.
If $n\leq 0$, we define $I_n(\varphi)=R$.
If $n\geq 1$, we define $I_n(\varphi)$ to be the image of 
$\varphi_n:\ext^nF\otimes\ext^n E^*\rightarrow R$ given by
\[
\varphi_n(f_1\wedge\cdots\wedge f_n\otimes \epsilon_1\wedge\cdots\wedge
\epsilon_n)=\det(\epsilon_i(\varphi(f_j))).
\]

Let $M$ be a finitely generated $R$-module.
Take a presentation
\begin{equation}\label{presentation.eq}
F\xrightarrow{\varphi} R^r\rightarrow M\rightarrow 0
\end{equation}
with $F$ being $R$-free (not necessarily finite).
For $j\in \Bbb Z$, $I_{r-j}(\varphi)$ is independent of the choice of 
the presentation (\ref{presentation.eq}).
We denote $I_{r-j}(\varphi)$ by $\Fitt_j(M)$, and call it the $j$th
Fitting ideal of $M$, see \cite[(20.2)]{Eisenbud}.

The construction of Fitting ideals commutes with base change
\cite[Corollary~20.5]{Eisenbud}.
So for a scheme $Y$ and a quasi-coherent $\O_Y$-module $\M$ of finite type
and $j\in \Bbb Z$, the quasi-coherent ideal $\uFitt_j(\M)$ of $\O_Y$ is defined
in an obvious way.
If $\M$ is locally of finite presentation, then $\uFitt_j(\M)$ is 
of finite type for any $j$.

\begin{lemma}\label{Fitting-base-change.thm}
Let $f:Y\rightarrow Z$ be a morphism of schemes, $\M$ a quasi-coherent
$\O_Z$-module of finite type, and $j\in \Bbb Z$.
Then $(\uFitt_j(\M))\O_Y=\uFitt_j(f^*(\M))$.
\end{lemma}

\proof Follows immediately by \cite[Corollary~20.5]{Eisenbud}.
\qed

As in \cite[Proposition~20.8]{Eisenbud}, we can prove the following easily.

\begin{lemma}\label{Fitting-locally-free.thm}
Let $Y$ be a scheme, $\M$ a quasi-coherent $\O_Y$-module of finite type,
and $j\geq 0$.
If $\uFitt_j(\M)=\O_Y$ and $\uFitt_{j-1}(\M)=0$, then 
$\M$ is locally free of well-defined rank $j$.
\end{lemma}

\paragraph Throughout the article, $S$ denotes a scheme, and $G$ denotes an
$S$-group scheme.
Throughout, $X$ denotes a $G$-scheme (i.e., an $S$-scheme with a left 
$G$-action).
We always 
assume that the second projection $p_2:G\times_S X\rightarrow X$ is
flat.

\paragraph
As in \cite[section~29]{ETI}, let $B_G^M(X)$ be the diagram
\[
\xymatrix{
G \times G\times X
\ar@<1.5em>[r]^-{1_G\times a}
\ar[r]^-{\mu\times 1_X}
\ar@<-1.5em>[r]^-{p_{23}} &
G\times X
\ar@<.75em>[r]^-a
\ar@<-.75em>[r]^-{p_2}
&
X
},
\]
where $a:G\times X\rightarrow X$ is the action, 
$\mu:G\times G\rightarrow G$ is the product, 
and $p_2:G\times X\rightarrow X$ and 
$p_{23}:G\times G\times X\rightarrow G\times X$ are
appropriate projections.
Note that $B_G^M(X)$ has flat arrows.
To see this, it suffices to see that $a:G\times X\rightarrow X$ is flat.
But since $a=p_2 b$ with $b$ an isomorphism, where $b(g,x)=(g,gx)$, this
is trivial.

By definition, a $(G,\Cal O_X)$-module is an $\Cal O_{B_G^M(X)}$-module.
It is said to be equivariant, quasi-coherent or coherent, if it is so as an
$\Cal O_{B_G^M(X)}$-module, see \cite{ETI}.
The category of equivariant (resp.\ quasi-coherent, 
coherent) $(G,\Cal O_X)$-modules is
denoted by $\EM(G,X)$ (resp.\ $\Qch(G,X)$, $\Coh(G,X)$).
In general, $\Qch(G,X)$ is closed under kernels, small 
colimits (in particular, 
cokernels), and extensions
in the category of $(G,\Cal O_X)$-modules.
In particular, it is an abelian category with the (AB5) condition
(\cite[Lemma~7.6]{ETI}).

\paragraph Let $X$ be as above.
We say that $(\Cal M,\Phi)$ is a $G$-linearized $\Cal O_X$-module
if $\Cal M$ is an $\Cal O_X$-module, and 
$\Phi:a^*\Cal M\rightarrow p_2^*\Cal M$ 
an isomorphism of $\Cal O_{G\times X}$-modules such that
\[
(\mu\times 1_X)^*\Phi: (\mu\times 1_X)^*a^*\M \rightarrow
(\mu\times 1_X)^*p_2^*\M
\]
agrees with the composite map
\begin{multline*}
(\mu\times 1_X)^*a^*\Cal M
\xrightarrow d
(1_G\times a)^*a^*\Cal M
\xrightarrow\Phi
(1_G\times a)^*p_2^*\Cal M\\
\xrightarrow d
p_{23}^*a^*\Cal M
\xrightarrow\Phi
p_{23}^*p_2^*\Cal M
\xrightarrow d
(\mu\times 1_X)^*p_2^*\Cal M,
\end{multline*}
where $d$'s are canonical isomorphisms.
We call $\Phi$ the $G$-linearization of $\Cal M$.

A morphism $\varphi:(\Cal M,\Phi)\rightarrow (\Cal N,\Psi)$ of 
$G$-linearized $\Cal O_X$-modules is an $\Cal O_X$-linear map 
$\varphi:\Cal M\rightarrow \Cal N$ such that 
$\Psi a^*\varphi=p_2^*\varphi\Phi$.
We denote the category of $G$-linearized $\Cal O_X$-modules by 
$\Lin(G,X)$.
The full subcategory of $G$-linearized quasi-coherent $\Cal O_X$-modules
is denoted by $\LQ(G,X)$.

For $\Cal M\in \EM(G,X)$, 
$(\Cal M_0,\alpha_{\delta_1(1)}^{-1}\alpha_{\delta_0(1)})$ is
in $\Lin(G,X)$, and this correspondence gives an equivalence.
With this equivalence, $\Qch(G,X)$ is equivalent to $\LQ(G,X)$.
See the proof of \cite[Lemma~9.4]{ETI}.

\paragraph\label{locally-noetherian.par}
If $X$ and $G\times X$ are quasi-compact quasi-separated, then 
$\Qch(G,X)$ is Grothendieck.
If, moreover, $X$ is noetherian, then $\Qch(G,X)$ is locally noetherian, 
and $\Cal M\in\Qch(G,X)$ is a noetherian object if and only if 
$\Cal M_0$ is coherent (\cite[Lemma~12.8]{ETI}).

\begin{example}
Let $k$ be a field, $G$ an affine algebraic group, and $A$ a $G$-algebra
(that is, a $k$-algebra on which $G$ acts).
We say that $M$ is a $(G,A)$-module if $M$ is a $G$-module, $M$ is an 
$A$-module, the $k$-space structures coming from the $G$-module 
structure and
the $A$-module structure agree, and the product $A\otimes M\rightarrow
M$ is $G$-linear.
It is known that the category of $(G,A)$-modules is equivalent to
the category of quasi-coherent $(G,\O_{\Spec A})$-modules.
\end{example}

\paragraph
The restriction functor $(?)_0:\Qch(G,X)\rightarrow \Qch(X)$ is faithful 
exact.
With this reason, we sometimes let $\Cal N\in \Qch(X)$ mean 
$\Cal M\in\Qch(G,X)$ if $\Cal N=\Cal M_0$.
For example, $\Cal O_X$ means the quasi-coherent 
$(G,\Cal O_X)$-module $\Cal O_{B_G^M(X)}$, since $(\Cal O_{B_G^M(X)})_0
=\Cal O_X$.
Note that $\Cal M\in\Qch(G,X)$ is coherent (i.e., it is in $\Coh(G,X)$) if
and only if $\Cal N=\Cal M_0\in\Coh(X)$, and no confusion will occur.

Let $\Cal M'$ be a subobject of $\Cal M\in\Qch(G,X)$.
Then $\Cal M'_0\subset \Cal M_0$ and the $(G,\Cal O_X)$-module structure
of $\Cal M$ together determine the $(G,\Cal O_X)$-submodule structure of 
$\Cal M'$ uniquely.
This is similar to the fact that for a ring $A$ and an $A$-module $M$ and
its $A$-submodule $N$, the subset $N\subset M$ and the $A$-module structure
of $M$ together determine the $A$-submodule structure of $N$ uniquely.
So by abuse of notation, we sometimes say that $\Cal M'_0$ is a 
quasi-coherent $(G,\Cal O_X)$-submodule of $\Cal M_0$, 
instead of saying that $\Cal M'$ is a quasi-coherent 
$(G,\Cal O_X)$-submodule of $\Cal M$.
Applying this abuse to $\Cal O_X$, we sometimes say that $\Cal I\subset
\Cal O_X$ is a quasi-coherent $G$-ideal of $\Cal O_X$ 
(i.e., a quasi-coherent $(G,\Cal O_X)$-submodule of $\Cal O_X$).

\paragraph Let $(\M,\Phi)$ be a $G$-linearized $\O_X$-module,
and $\N$ an $\O_X$-submodule.
We identify $p_2^*\N$ by its image in $p_2^*\M$, since $p_2$ is flat.
Similarly, $a$ is also flat, and we identify $a^*\N$ by its image in $a^*\M$.
Then $\N$ is a $(G,\O_X)$-submodule if and only if 
$\Phi(a^*\N)=p_2^*\N$,
since then $\Phi:a^*\N\rightarrow p_2^*\N$ is an isomorphism.

So for a $G$-equivariant $\O_X$-module $\M$ and an $\O_X$-submodule $\N$ of
$\M_0$, 
$\N$ is a $(G,\O_X)$-submodule if and only if the image of
$a^*\N$ by the map $\alpha:a^*\M_0\rightarrow \M_1$
agrees with the image of $p_2^*\N$ by the map $\alpha:p_2^*\M_0\rightarrow 
\M_1$.

\begin{lemma}\label{Fitting-G-ideal.thm}
Let $\M$ be a quasi-coherent $(G,\O_X)$-module of finite type, 
and $j\in\Bbb Z$.
Then the Fitting ideal $\uFitt_j\M$ is a $G$-ideal of $\O_X$.
\end{lemma}

\proof By Lemma~\ref{Fitting-base-change.thm},
the two 
extended ideals $(\uFitt_j\M)\O_{G\times X}$ via $a:G\times X\rightarrow X$
and $p_2:G\times X\rightarrow X$ agree, since the former one is
$\uFitt_j(a^*\M)$, the latter one is $\uFitt_j(p_2^*\M)$, and
$a^*\M\cong p_2^*\M$.
\qed

\paragraph
The restriction $(?)_0:\Qch(G,X)\rightarrow\Qch(X)$ 
has a right adjoint, if the second projection 
$p_2:G\times X\rightarrow X$ is quasi-compact quasi-separated
(\cite[Lemma~12.11]{ETI}).

\section{Primary decompositions over noetherian schemes}

Let $Y$ be a scheme.
An ideal of $\Cal O_Y$ means a quasi-coherent ideal sheaves, unless otherwise
specified.
An $\Cal O_Y$-module means a quasi-coherent
module, unless otherwise specified.

\paragraph An ideal $\Cal P$ of $\Cal O_Y$ is said to be prime if 
$V(\Cal P)$ is integral.
An ideal $\Cal P$ of $\Cal O_Y$ is said to be quasi-prime if
$\Cal P\neq \Cal O_Y$, and if $\I$ and $\J$ are ideals of $\Cal O_Y$ such
that $\I\J\subset\Cal P$, then $\I\subset\Cal P$ or $\J\subset \Cal P$ 
holds.

\begin{lemma} 
Let $\Cal P$ be an ideal of $\Cal O_Y$.
If $\Cal P$ is prime, then $\Cal P$ is quasi-prime.
If $Y$ is quasi-compact quasi-separated and $\Cal P$ is quasi-prime,
then $\Cal P$ is prime.
\end{lemma}

\proof 
Replacing $Y$ by $V(\Cal P)$, $\Cal P$ by $0$, $\Cal I$ by $\Cal I\Cal O_{V
(\Cal P)}$, and $\Cal J$ by $\Cal J\Cal O_{V(\Cal P)}$, 
we may assume that $\Cal P=0$.

We prove the first part.
Since $Y$ is integral, it is irreducible and hence is non-empty.
Thus $\Cal O_Y\neq 0=\Cal P$.

Let $\eta$ be the generic point of $Y$.
Since $\Cal I\Cal J=0$, 
$\Cal I_\eta\Cal J_\eta=0$.
Since $\Cal O_{Y,\eta}$ is an integral domain, 
$\Cal I_\eta=0$ or $\Cal J_\eta=0$.
If $\Cal I_\eta=0$, then
\[
\Cal I\subset\Cal I_\eta\cap \Cal O_Y=0\cap \Cal O_Y.
\]
This is zero by (\ref{reduced-sim.par}), applied to $\Spec \O_{Y,\eta}
\rightarrow Y$.
Similarly, $\Cal J_\eta=0$ implies $\Cal J=0$.

We prove the second part.
Since $Y$ is quasi-compact, it has a finite affine open covering 
$(U_i)_{i=1}^n$.
We may 
assume that $U_1,\ldots,U_s$ are reduced, and $U_i$ is not reduced for
$i>s$.
Let $\I_i=0_{U_i}\cap \Cal O_Y$ be the pull-back of zero for $i\leq s$.
For $i>s$, there is a non-zero ideal $\J_i$ of $\Cal O_{U_i}$ such that
$\J_i^2=0$, since $U_i$ is affine and non-reduced.
Set $\I_i=(\J_i\cap \Cal O_Y)^2$.

Since the inclusion $U_i\hookrightarrow Y$ is quasi-compact,
$\I_i|_{U_i}=0$.
Thus $\I_1\cdots\I_n=0$.
By assumption, there exists some $i$ such that $\I_i=0$.
If $i>s$, then $\J_i\cap \Cal O_Y=0$ by assumption again.
So $\J_i=(\J_i\cap \Cal O_Y)|_{U_i}=0$, and this is a contradiction.
So $i\leq s$.
Thus $Y$ is the scheme theoretic image of the inclusion $U_i\hookrightarrow
Y$.
By (\ref{reduced-sim.par}), $Y$ is reduced.

Since $\Cal O_Y\neq 0$, $Y$ is non-empty.
Assume that $Y$ is not irreducible.
Then $Y=Y_1\cup Y_2$ for some closed subsets $Y_i\neq Y$.
Let us consider the reduced structure of $Y_i$, and set $\K_i=\Cal I(Y_i)$.
Then $\K_1\cap \K_2=0$.
By assumption, $\K_1=0$ or $\K_2=0$.
This contradicts $Y_i\neq Y$.
So $Y$ must be irreducible, and $Y$ is integral.
\qed

\paragraph \label{maximal.par}
An ideal $\Cal M$ of $\Cal O_Y$ is said to be maximal, if
$\Cal M$ is a maximal element of the set of proper (i.e., not equal to 
$\Cal O_Y$) ideals of $\Cal O_Y$.
The defining ideal of a closed point (with the reduced structure) is maximal.
Since a non-empty quasi-compact $T_0$-space has a closed point, 
if $Y$ is non-empty and quasi-compact, then $\Cal O_Y$ has a maximal ideal.
It is easy to see that a maximal ideal is a prime ideal.

\begin{lemma}\label{prime-map.thm}
If $f:Y\rightarrow Z$ is a morphism and $\Cal P$ is a prime
ideal of $\Cal O_Y$, then $\Cal P\cap\Cal O_Z$ is a prime ideal.
\end{lemma}

\proof This is nothing but the restatement of the fact that
the scheme theoretic image of the composite
\[
V(\Cal P)\hookrightarrow Y\xrightarrow{f}Z
\]
is integral, see (\ref{reduced-sim.par}).
\qed

\paragraph Let $\I$ be an ideal of $\Cal O_Y$.
For an affine open subset $U$ of $Y$, we define $\Gamma(U,\sqrt{\I}):=
\sqrt{\Gamma(U,\I)}$.
This defines a quasi-coherent ideal $\sqrt{\I}$ of $\O_Y$.
We call $\sqrt{\I}$ the {\em radical} of $\I$.
The formation of $\sqrt{\I}$ is local (that is, for any open subset $U$ of
$Y$, $\sqrt{\I}|_U=\sqrt{\I|_U}$), and $V(\I)$ is reduced if and only
if $\I=\sqrt{\I}$.
Hence $\Cal P=\sqrt{\Cal P}$ for a prime ideal $\Cal P$ of $\Cal O_Y$.
Note also that $\sqrt{\I}=\O_X$ if and only if $\I=\O_X$,
since the formation of $\sqrt{\I}$ is local.

\begin{lemma}
For an ideal $\I$ of $\Cal O_Y$, 
\[
\sqrt{\I}=(\bigcap \Cal P)^\star,
\]
where the intersection is taken over all prime ideals $\P$ containing $\I$.
\end{lemma}

\proof Note that the assertion is well-known for affine schemes.
Set the right hand side to be $\J$.
Since $\I\subset\cP$, we have $\sqrt{\I}\subset\sqrt{\cP}=\cP$.
So $\sqrt{\I}\subset\J$ is obvious.

To prove $\sqrt{\I}\supset \J$, it suffices to prove that
$\sqrt{\I}|_U\supset \J|_U$ for any affine open subset $U$ of $Y$.
Let $\Q$ be a prime ideal of $\Cal O_U$ such that $\Q\supset \I|_U$.
Then $\Q\cap \Cal O_Y$ is prime by Lemma~\ref{prime-map.thm},
and $\Q\cap \Cal O_Y\supset \I|_U\cap \Cal O_Y\supset \I$.
Hence $\Q\cap \Cal O_Y\supset \J$.
Thus $\Q\supset (\Q\cap\Cal O_Y)|_U\supset \J|_U$.
Hence $\J|_U\subset \bigcap \Q=\sqrt{\I|_U}=\sqrt{\I}|_U$.
\qed

\paragraph An ideal $\I$ of $\Cal O_Y$ is said to be a {\em radical ideal}
if $\sqrt{\I}=\I$.

\begin{lemma}\label{radical-intersect.thm}
Let $(\I_\lambda)_{\lambda\in\Lambda}$ be a family of radical ideals of 
$\Cal O_Y$.
Then $(\bigcap_\lambda \I_\lambda)^\star$ is a radical ideal.
\end{lemma}

\proof 
Let $Y_\lambda:=V(\I_\lambda)$, and $Y':=\coprod_\lambda Y_\lambda$.
Note that $Y'$ is reduced.
So the scheme theoretic image of the canonical map $\varphi:Y'\rightarrow Y$
is also reduced by (\ref{reduced-sim.par}).
On the other hand, the defining ideal of $\Sim(\varphi)$ is
$(\bigcap_\lambda \I_\lambda)^\star$, and we are done.
\qed

\paragraph \label{def-ass.par}
From now, until the end of this section, let $Y$ be noetherian.
For an $\Cal O_Y$-module $\Cal M$, the subset 
\[
\{y\in Y\mid \Hom_{\Cal O_{Y,y}}(\kappa(y),\Cal M_y)\neq 0\}
\]
of $Y$ is denoted by $\Ass(\Cal M)$.
A point $y$ in $\Ass(\Cal M)$ is called an associated point of $\Cal M$.
The closed subscheme $V(\overline{\{y\}})$ with $y\in\Ass(\Cal M)$ is 
called an associated component.
The defining ideal of an associated component is called an associated prime
ideal.

\begin{lemma}
Let $\M$ be a coherent $\Cal O_Y$-module, and $\N$ its submodule.
Then the following are equivalent.
\begin{description}
\item[(i)] $\Ass(\M/\N)=\{y\}$ is a singleton.
\item[(ii)] $\M\neq\N$, and
if $\cL$ is a submodule of $\M$, $\I$ an ideal of $\Cal O_Y$, 
$\I\cL\subset\N$, and $\cL\not\subset\N$, then $\I\subset\sqrt{\N:\M}$.
\end{description}
If this is the case, $\sqrt{\N:\M}$ is a prime ideal, 
and $y$ is the generic point of $V(\sqrt{\N:\M})$.
\end{lemma}

\proof Left to the reader.
\qed

\paragraph
Let $\Cal M$ be a coherent $\Cal O_Y$-module and $\Cal N$ its submodule.
If the equivalent conditions of the lemma is satisfied, then we say that
$\Cal N$ is a {\em primary submodule} of $\Cal M$.
If $\Cal M=\O_Y$, then we say that $\Cal N$ is a {\em primary ideal}.
If $0$ is a primary submodule of $\Cal M$, then we say that $\Cal M$ is
a {\em primary module}.
If $\Cal O_Y$ is a primary module, then we say that $Y$ is {\em primary}.

If $\Cal N$ is a primary submodule of $\Cal M$, $\Cal M$ is coherent, and
$\Ass(\M/\N)=\{y\}$, then we say that $\N$ is $y$-primary.
We also say that $\N$ is $\sqrt{\N:\M}$-primary.

Note that if $\N$ is a primary submodule of $\M$, then $\N:\M$ is a
primary ideal (easy).

We consider that $Y$ is an ordered set with respect to the order given by
$y\leq y'$ if and only if $y$ is a generalization of $y'$.
For a coherent $\Cal O_Y$-module $\Cal M$, 
any minimal element of $\Supp \Cal M$ is a member of $\Ass (\Cal M)$.
The set of minimal elements of $\Supp\Cal M$ is denoted by 
$\Min(\Cal M)$.
If $\Cal M$ is coherent, then $\Ass(\Cal M)$ is a finite set.
If $\Cal M$ is quasi-coherent, then $\Ass(\Cal M)=\emptyset$ if and only if
$\Cal M=0$.
If $\Cal M$ is quasi-coherent, then $\Supp(\M)\supset 
\Ass(\Cal M)\supset\Min(\M)$.

\paragraph
There is a one-to-one correspondence between
coherent prime ideals of $\Cal O_Y$ 
and integral closed subschemes of $Y$ ($\Cal P$ corresponds to
$V(\Cal P)$).
So they are also in one to one correspondences with points in $Y$
($V(\Cal P)$ corresponds to its generic point).
So for a coherent $\Cal O_Y$-module $\Cal M$, $\Supp \M$, $\Ass\Cal M$, and
$\Min \Cal M$ are sometimes considered as sets of prime ideals of $Y$.
An element of $\Ass\Cal M\setminus\Min\Cal M$ is called an embedded prime
ideal of $\Cal M$.

\paragraph \label{primary.par}
Let $\Cal M$ be a coherent sheaf over $Y$,
$\Ass(\Cal M)=\{y_1,\ldots,y_r\}$, and $y_1,\ldots,y_r$ be distinct.
Let $0= M_{i,1}\cap\cdots \cap M_{i,r_i}$ 
be a minimal primary decomposition of 
$0\subset\Cal M_{y_i}$.
Since $\depth\Cal M_{y_i}=0$, by reordering if necessary,
we may assume that $M_{i,1}$ is $\frak m_{y_i}$-primary.
We set $L_i=M_{i,1}$.
Note that $H^0_{\frak m_{y_i}}(\Cal M_{y_i})=\bigcap_{j\geq2}M_{i,j}$.
Hence $H^0_{\frak m_{y_i}}(\Cal M_{y_i})\cap L_i=0$ and
$\Cal M_{y_i}/L_i$ has a finite length.

Let $\varphi_i:\Spec\Cal O_{y_i}\rightarrow Y$ be the canonical map.
Let $\Cal N_i$ be the kernel of the composite
\[
\Cal M\rightarrow (\varphi_i)_*\Cal M_{y_i}\rightarrow
(\varphi_i)_*(\Cal M_{y_i}/L_i)
\]
so that there  is a monomorphism $\Cal M/\Cal N_i\hookrightarrow
(\varphi_i)_*(\Cal M_{y_i}/L_i)$.
It is easy to see that $\Ass(\Cal M/\Cal N_i)=\{y_i\}$.

Since $(\Cal N_i)_{y_i}= L_i$ and $H^0_{\frak m_{y_i}}(L_i)=0$, 
we have that $y_i\notin \Ass(\Cal N_i)$.
In particular, $\Ass(\Cal N_1\cap \cdots\cap \Cal N_r)=\emptyset$.

So there is a decomposition
\[
0=\Cal N_1\cap\cdots\cap \Cal N_r
\]
such that $\Ass(\Cal M/\Cal N_i)=\{y_i\}$.
We call such a decomposition a minimal primary decomposition of $0$.

If $\Cal N$ is a coherent $\Cal O_Y$-submodule of $\Cal M$ and 
$\Ass(\Cal M/\Cal N)=\{y_1,\ldots,y_r\}$
($y_1,\ldots,y_r$ are distinct),
then
there is a decomposition
\begin{equation}\label{N-M.eq}
\Cal N=\Cal M_1\cap\cdots\cap \Cal M_r
\end{equation}
such that $\Ass(\Cal M/\Cal M_i)=\{y_i\}$.
We call such a decomposition a minimal primary decomposition of $\Cal N$.
Such a decomposition is not unique in general.
If a coherent $\Cal O_X$-submodule $\Cal W$ agrees with some $\Cal M_i$ for
some minimal primary decomposition (\ref{N-M.eq}), we say that
$\Cal W$ is a primary component of $\Cal N$.
Note that a primary component $\M_i$ for $y_i\in \Min(\M/\N)$ is known
to be unique.
If $\M=\O_X$ and $\N$ is a radical ideal,
then $\O_X/\N$ does not have an embedded prime ideal, and $\M_i$ is the unique
prime ideal such that the generic point of $V(\M_i)$ is $y_i$.

\begin{lemma}\label{M-ann.thm}
Let $\Cal M$ be a coherent $\Cal O_Y$-module, and $\Cal N$ a coherent 
submodule of $\Cal M$.
Let {\rm(\ref{N-M.eq})} be a minimal primary decomposition of $\Cal N$.
If $\Cal M/\Cal N$ does not have an embedded prime ideal, then
\begin{equation}\label{N-M2.eq}
\Cal N:\Cal M=\bigcap_{i=1}^r \Cal M_i:\Cal M
\end{equation}
is a minimal primary decomposition.
In particular, $\Ass \Cal M/\Cal N=\Ass\Cal O_X/(\Cal N:\Cal M)$.
\end{lemma}

\proof It is obvious that the equation (\ref{N-M2.eq}) 
holds, and it is a primary 
decomposition.

Since (\ref{N-M.eq}) is minimal, $\sqrt{\Cal M_i:\Cal M}$ are distinct, and
there is no incidence relation each other.
The minimality of (\ref{N-M2.eq}) follows easily.
\qed

\paragraph
Let $\Cal M$ be a coherent sheaf over $Y$.
Note that $\Cal M$ satisfies Serre's $(S_1)$-condition 
(see \cite[(5.7.2)]{EGA-IV}) if and only if $\Cal M$ has no embedded prime 
ideal.
So $\M$ is primary if and only if $\M$ satisfies the $(S_1)$ and
$\Supp\M$ is irreducible.

\begin{lemma}\label{S_1-flat.thm}
Let $U$ be an open subset of $Y$, and $i:U\hookrightarrow Y$ be the
inclusion.
Let $\M$ be a coherent sheaf over $Y$.
\begin{description}
\item[(i)] $\Ass(i^*\M)=\Ass(\M)\cap U$.
\item[(ii)] $\Ass(\uH^0_{Y\setminus U}(\M))=\Ass(\M)\setminus U$,
where $\uH^0_{Y\setminus U}(\M)$ is the kernel of 
the canonical map $\M\rightarrow i_*i^*\M$.
\item[(iii)] The following are equivalent:
\begin{description}
\item[(a)] $\Ass(\M)\subset U$;
\item[(b)] $\Ass(\M)=\Ass(i^*\M)$;
\item[(c)] $\uH^0_{Y\setminus U}(\M)=0$;
\item[(d)] $\M\rightarrow i_*i^*\M$ is monic,
\end{description}
\item[(iv)] Assume that $U\cap \Supp\M$ is dense in $\Supp M$.
If $\M$ satisfies $(S_1)$, then $\M\rightarrow i_*i^*\M$ is monic.
\end{description}
\end{lemma}

\proof {\bf (i)} This is because $(i^*\M)_x\cong \M_x$ for $x\in U$.

{\bf (ii)} Note that $\uH^0_{Y\setminus U}\Cal M$ is coherent.
Let $x\in Y\setminus U$.
The image of the injective map $(\uH^0_{Y\setminus U}\Cal M)_x
\rightarrow \M_x$ is identified with $H^0_{\I_x}(\M_x)$, where $\I$ 
is any ideal such that $V(\I)=Y\setminus U$ set theoretically.
As any map $\kappa(x)\rightarrow \M_x$ factors through $H^0_{\I_x}(\M_x)$,
the assertion follows.

{\bf (iii)} {\bf (a)$\Leftrightarrow$(b)} follows from {\bf (i)}.
{\bf (c)} is equivalent to $\Ass(\uH^0_{Y\setminus U}(\M))=\emptyset$.
By {\bf (ii)}, {\bf (a)$\Leftrightarrow$(c)} follows.
{\bf (c)$\Leftrightarrow$(d)} is trivial.

{\bf (iv)}
Note that any associated point of $\M$ in $Y\setminus U$
is embedded by assumption.
As $\M$ does not have an embedded prime ideal by assumption, 
$\Ass\M\subset U$.
By {\bf (iii)}, $\M\rightarrow i_*i^*\M$ is monic.
\qed

\begin{corollary}\label{S_1.thm}
Let 
$f:Z\rightarrow Y$ be a flat morphism of finite type.
If the image of $f$ is dense in $Y$ and $Y$ satisfies $(S_1)$, then
$\Sim f=Y$.
\end{corollary}

\proof Replacing $Z$ by $\Im f$, we may assume that $f$ is an open immersion.
The assertion follows immediately by the lemma applied to $\Cal O_Y$.
\qed

\begin{lemma}\label{primary-sim.thm}
Let $f:Z\rightarrow Y$ be a morphism of 
noetherian schemes.
Let $\M$ be a primary coherent $\O_Z$-module.
If $\N$ is a coherent $\O_Y$-submodule of $f_*\M$, then $\N$ is
either zero or primary.
In particular, if $Z$ is primary, then $\Sim f$ is primary.
\end{lemma}

\proof 
First, note that if $A\rightarrow B$ is an injective homomorphism
between noetherian rings, $M$ a primary finitely generated faithful 
$B$-module, and $N$ a finitely generated
$A$-submodule of $M$, then either $N=0$ or
$N$ is primary and $\Supp N=\Spec A$, set theoretically.

Indeed, if $n\in N\setminus 0$, then $0:_B n\subset \sqrt{0_B}$.
So $0:_A n\subset \sqrt{0_A}$.

To prove the lemma, replacing $Z$ by $V(\uann\M)$, we may assume that
$\uann\M=0$.
By Lemma~\ref{M-ann.thm}, $Z$ is primary.
In particular, $Z$ is irreducible.
Next, replacing $Y$ by $\Sim f$, we may assume that $Y=\Sim f$.
In particular, $Y$ is irreducible, and $\O_Y\rightarrow f_*\O_Z$ is monic.
We may assume that $\N\neq 0$.
Let $U=\Spec A$ be an open subset of $Y$ such that $i^*\N\neq 0$,
where $i:U\hookrightarrow Y$ is the inclusion.
Let $V=\Spec B$ be a non-empty open subset of $f^{-1}(U)$.
By Lemma~\ref{S_1-flat.thm}, 
\[
i^*\N\rightarrow i^*f_*\M\cong g_*j^*\M\rightarrow 
g_*k_*k^*j^*\M=(gk)_*(jk)^*\M
\]
is monic, where $j:f^{-1}(U)\rightarrow Z$ is the
inclusion, $k:V\rightarrow f^{-1}(U)$ is the inclusion, 
and $g:f^{-1}(U)\rightarrow U$ is the 
restriction of $f$.
By the affine case above, $\Ass(i^*\N)=\{\eta\}$, where $\eta$ is the
generic point of $Y$.

The diagram
\[
\xymatrix{
\N \ar[d] \ar[r] & i_*i^*\N \ar[d] \\
f_*\M \ar[r] & f_*(jk)_*(jk)^*\M=i_*(gk)_*(jk)^*\M
}
\]
is commutative.
Thus the canonical map $\N\rightarrow i_*i^*\N$ is monic by 
Lemma~\ref{S_1-flat.thm} applied to $\M$.
By Lemma~\ref{S_1-flat.thm} applied to $\N$, 
$\Ass\N=\Ass(i^*\N)=\{\eta\}$, and $\N$ is primary, as desired.
\qed

\section{$G$-prime and $G$-radical $G$-ideals}

Let $S$ be a scheme, $G$ an $S$-group scheme, and
$X$ a $G$-scheme.
We say that $X$ is a p-flat $G$-scheme if 
the second projection 
$p_2:G\times X \rightarrow X$ is flat.
If $G$ is flat over $S$, then any $G$-scheme is p-flat.
We always assume that $X$ is p-flat.
Although we do not assume that $G$ is $S$-flat, the sheaf theory as in 
\cite{ETI} and \cite{LCDS} goes well, since we assume that $X$ is p-flat
and hence $B_G^M(X)$ has flat arrows.

In the rest of the paper, an $\Cal O_X$-module and an ideal of $\Cal O_X$
are required to be quasi-coherent, unless otherwise specified.
A $(G,\Cal O_X)$-module and a $G$-ideal of $\Cal O_X$ are
also required to be quasi-coherent unless otherwise specified.

\paragraph Let $\Cal M$ be a $(G,\Cal O_X)$-module.
Note that the sum $\sum_\lambda\Cal M_\lambda$ of quasi-coherent 
$(G,\Cal O_X)$-submodules
$\Cal M_\lambda$ is a quasi-coherent $(G,\Cal O_X)$-submodule.
If $\Cal N$ and $\Cal L$ are quasi-coherent $(G,\Cal O_X)$-submodules,
then $\Cal N\cap \Cal L$ is again a quasi-coherent $(G,\Cal O_X)$-submodule.
For a quasi-coherent $G$-ideal $\Cal I$, $\Cal I\Cal N$ is 
a quasi-coherent $(G,\Cal O_X)$-submodule.
If, moreover, $\Cal I$ is coherent, then being the kernel of the canonical map
\begin{equation}\label{colon.eq}
\M\rightarrow\uHom_{\Cal O_X}(\I,\M/\N),
\end{equation}
$\Cal N:\Cal I$
is also a quasi-coherent $(G,\Cal O_X)$-submodule, see
\cite[(7.11)]{ETI} and \cite[(7.6)]{ETI}.
More generally,

\begin{lemma}
Let $Y$ be a scheme, $\M$ a quasi-coherent $\O_Y$-module, 
$\N$ a quasi-coherent $\O_Y$-submodule of $\M$, and 
$\I$ a quasi-coherent ideal of $\O_Y$.
If $\I$ is of finite type, then $\N:\I$, the kernel of {\rm(\ref{colon.eq})},
is a quasi-coherent submodule of $\M$.
If $X$ is a $G$-scheme, $\M$ a quasi-coherent $(G,\O_X)$-module,
$\N$ its quasi-coherent $(G,\O_X)$-submodule, 
and $\I$ a quasi-coherent $G$-ideal
of finite type, 
then $\N:\I$ is a
quasi-coherent $(G,\O_X)$-submodule of $\M$.
\end{lemma}

\proof We prove the first assertion.
For an affine open subset $U$ of $Y$, 
$\Gamma(U,\N:\I)$ is the kernel of $M\rightarrow \Hom_A(I,M/N)$, where
$A:=\Gamma(U,\O_X)$, $I:=\Gamma(U,\I)$, $M:=\Gamma(U,\M)$, and
$N:=\Gamma(U,\N)$.
So $\Gamma(U,\N:\I)=N:I$.
Since $(N:I)B=N\otimes_AB:IB$ for a flat $A$-algebra $B$, the formation
of a colon module (for of finite type $\I$) 
is compatible with the localization.
So $\N:\I$ is quasi-coherent.

Next, we prove the second assertion.
By the reason above, formation of a colon module (for of finite type $\I$)
is compatible with a flat base change.
So
\[
\Phi(a^*(\N:\I))=\Phi(a^*\N:a^*\I)=\Phi(a^*\N):a^*\I=p_2^*\N:p_2^*\I
=p_2^*(\N:\I).
\]
This shows that $\N:\I$ is a $(G,\O_X)$-submodule of $\M$.
\qed

Similarly, if $\Cal L$ is a quasi-coherent $(G,\O_X)$-submodule of
$\M$ of finite type, then $\Cal N:\Cal L$ is a quasi-coherent
$G$-ideal of $\Cal O_X$.

For an $\Cal O_X$-submodule $\frak m$, quasi-coherent or not, of $\Cal M$,
the sum of all 
quasi-coherent $(G,\Cal O_X)
$-submodules of $\Cal M$ contained in $\frak m$ is the
largest quasi-coherent $(G,\Cal O_X)$-submodule of $\Cal M$ 
contained in $\frak m$.
We denote this by $\frak m^*$.
If $\frak a$ is a quasi-coherent ideal of $\Cal O_X$ 
and $Y=V(\frak a)$, then
we denote $V(\frak a^*)$ by $Y^*$.
$Y^*$ is the smallest closed $G$-subscheme of $X$ containing $Y$.

For a morphism $f:Y \rightarrow X$,
$\Ker(\Cal O_X\rightarrow f_*\Cal O_Y)^*$
defines the smallest closed $G$-subscheme $Y'$ of $X$ such that
$f^{-1}(Y')=Y$.
We call $Y'$ the $G$-scheme theoretic image of $Y$ by $f$, and 
denote it by $\GSim(f)$.
Clearly, $\GSim(f)\supset\Sim(f)$ and $\GSim(f)=\Sim(f)^*$.
Note that for a closed subscheme $Y$ of $X$, $Y^*$ is the
$G$-scheme theoretic image of the inclusion $Y\hookrightarrow X$.
It is easy to verify that, for a closed subscheme $Y$ of $X$, 
the $G$-scheme theoretic image of the action
$G\times Y \rightarrow X$ ($(g,y)\mapsto gy$) is $Y^*$.

If $f:Y\rightarrow X$ is a quasi-compact quasi-separated $G$-morphism of
$G$-schemes, 
then $\Ker(\Cal O_X\rightarrow f_*\Cal O_Y)$ is a quasi-coherent $G$-ideal.
So $\GSim(f)=\Sim(f)=V(\Ker(\Cal O_X\rightarrow f_*\Cal O_Y))$.

\begin{lemma}\label{GSim.thm}
Let $f:V\rightarrow X$ be a $G$-morphism of $G$-schemes.
Let $Y$ be a closed subscheme of $V$.
Then
\[
\GSim(Y^*\hookrightarrow V\rightarrow X)=
\Sim(Y\hookrightarrow V\rightarrow X)^*.
\]
\end{lemma}

\proof
\[
\GSim(Y^*\hookrightarrow V\rightarrow X)\supset
\GSim(Y\hookrightarrow V\rightarrow X)
=
\Sim(Y\hookrightarrow V\rightarrow X)^*.
\]

We prove the opposite inclusion.
$f^{-1}(\Sim(Y\hookrightarrow V\rightarrow X)^*)$ is a $G$-closed
subscheme of $V$ containing $Y$.
So it also contains $Y^*$ by the minimality of $Y^*$.
By the minimality of
$\GSim(Y^*\hookrightarrow V\rightarrow X)$, 
we have 
$\GSim(Y^*\hookrightarrow V\rightarrow X)\subset
\Sim(Y\hookrightarrow V\rightarrow X)^*$.
\qed

\begin{lemma}\label{intersection.thm}
Let $(\frak m_\lambda)_{\lambda\in\Lambda}$ be a family of
$\Cal O_X$-submodules of $\Cal M$.
Then $(\bigcap_{\lambda}\fm_\lambda^*)^*=
(\bigcap_{\lambda}\fm_\lambda)^*$.
\end{lemma}

\proof Since $\fm_\lambda^*\subset\fm_\lambda$ for each $\lambda$,
we have $(\bigcap_{\lambda}\fm_\lambda^*)^*\subset
(\bigcap_{\lambda}\fm_\lambda)^*$.

On the other hand, since $\fm_\lambda^*\supset(\bigcap_\lambda\fm_\lambda)^*$
for each $\lambda$, we have
$\bigcap\fm_\lambda^*\supset(\bigcap_\lambda\fm_\lambda)^*$.
By the maximality of 
$(\bigcap\fm_\lambda^*)^*$, 
we have
$(\bigcap\fm_\lambda^*)^*\supset
(\bigcap_\lambda\fm_\lambda)^*$.
\qed

\begin{corollary}\label{intersection2.thm}
Let $\frak m$ and $\frak n$ be $\Cal O_X$-submodules of $\Cal M$.
Then $\frak m^*\cap\frak n^*=(\frak m\cap\frak n)^*$.
\end{corollary}

\proof
Follows immediately from the lemma, since
$\fm^*\cap\fn^*=(\fm^*\cap\fn^*)^*$.
\qed

\begin{lemma}\label{colon.thm}
Let $\frak m$ be an $\Cal O_X$-submodule of $\Cal M$.
If $\M$ is of finite type,
then $(\frak m:\Cal M)^*=\frak m^*:\Cal M$.
\end{lemma}

\proof Set $\Cal I:=(\frak m:\Cal M)^*$.
Then $\Cal I\Cal M\subset\frak m$ and $\Cal I\Cal M$ is a quasi-coherent
$(G,\Cal O_X)$-submodule of $\Cal M$.
Hence $\Cal I\Cal M\subset\frak m^*$ by the maximality,
and $\Cal I\subset \frak m^*:\Cal M$.

On the other hand, $\frak m^*:\Cal M\subset\frak m:\Cal M$.
Hence $\frak m^*:\Cal M\subset\Cal I$ by the maximality.
\qed

\begin{lemma}
Let $\frak m$ be an $\Cal O_X$-submodule of $\Cal M$,
and $\Cal J$ a finite-type $G$-ideal of $\Cal O_X$.
Then $(\frak m:\J)^*=\frak m^*:\J$.
\end{lemma}

\proof Similar.
\qed

\paragraph\label{ukya.par}
We denote the scheme $X$ with the trivial $G$ action by $X'$.
Thus $G\times X'$ is the principal $G$-bundle (i.e., the $G$-scheme with
the $G$-action given by $g(g',x)=(gg',x)$).

\paragraph Let us consider the diagram
\[
\xymatrix{
G\times G\times G\times X~
\ar@<1.5em>[r]^-{1_G\times \mu\times 1_X}
\ar[r]^-{\mu\times 1_{G\times X}}
\ar@<-1.5em>[r]^-{p_{234}}
&
~G\times G\times X
\ar@<.75em>[r]^-{\mu\times 1}
\ar@<-.75em>[r]^-{p_{23}}
& 
G\times X
\ar[r]^-{p_2}
&
X
}
\]
on the finite category $\Delta^+_M$ (see for the definition, 
\cite[(9.1)]{ETI}).
For $\Cal M\in\Qch(X)$, $\Bbb A\Cal M$ is in $\Qch(G,X)$, where 
$\Bbb A=(?)_{\Delta_M}\circ L_{[-1]}$ is the ascent functor 
\cite[(12.9)]{ETI}.
Thus we may say that $p_2^*\Cal M$ is a quasi-coherent $(G,\Cal O_X)$-module,
since $(\Bbb A\Cal M)_0=p_2^*\Cal M$.
The $G$-linearization of $p_2^*\Cal M$ is the canonical isomorphism
$d:(\mu\times 1)^*p_2^*\Cal M\rightarrow p_{23}^*p_2^*\Cal M$.

\paragraph Let $a:G\times X'\rightarrow X$ be the action,
where $X'$ is as in (\ref{ukya.par}).
Then $a$ is a $G$-morphism.
Thus $a^*\Cal M$ is a quasi-coherent 
$G$-linearized $\Cal O_X$-module for $\Cal M\in\Qch(G,X)$.
The $G$-linearization is the composite map
\[
(\mu\times 1)^*a^*\Cal M
\xrightarrow d
(1\times a)^*a^*\Cal M
\xrightarrow\Phi
(1\times a)^*p_2^*\Cal M
\xrightarrow d
p_{23}^*a^*\Cal M.
\]
Since $\Phi:a^*\Cal M\rightarrow p_2^*\Cal M$ is a $G$-linearization,
$\Phi:a^*\Cal M\rightarrow p_2^*\Cal M$ is an isomorphism of $G$-linearized
$\Cal O_X$-modules.
In particular, the composite map
\[
\omega: \Cal M\xrightarrow u a_*a^*\Cal M\xrightarrow \Phi a_*p_2^*\Cal M
\]
is $(G,\Cal O_X)$-linear.

\begin{example}
Consider the case that everything is affine:
$S=\Spec R$, $G=\Spec H$, and $X=\Spec A$.
Let $\M$ be a quasi-coherent $(G,\Cal O_X)$-module.
Then $M=\Gamma(X,\M)$ is a $(G,A)$-module.
$p_2^*\Cal M$ corresponds to the $H\otimes A$-module $H\otimes M$
(the action of $H\otimes A$ on $H\otimes M$ is given by
$(h\otimes a)(h'\otimes m)=hh'\otimes am$).
$a^*\Cal M$ corresponds to the $H\otimes A$-module $H\otimes M$
(the action of $H\otimes A$ on $H\otimes M$ is given by
$(h\otimes a)(h'\otimes m)=\sum_{(a)}hh'a_{(1)}\otimes a_{(0)}m$,
where we employ Sweedler's notation, that is, $\omega_A(a)=\sum_{(a)}a_{
(0)}\otimes a_{(1)}$, where $\omega_A:A\rightarrow A\otimes H$ is the
coaction of the right $H$-comodule $A$ corresponding to the $G$-module $A$.
Then $\Phi:a^*\M\rightarrow p_2^*\M$ corresponds to the map
$\Box:H\otimes M\rightarrow H\otimes M$ given by $\Box(h\otimes m)
=\sum_{(m)}h(\Cal Sm_{(1)})\otimes m_{(0)}$, where $\Cal S$ is the antipode
of the Hopf algebra $H$.
Thus $\omega:\M\rightarrow a_*p_2^*\M$ corresponds to the map
$\omega':M\rightarrow H\otimes M$ given by $\omega'(m)=\sum_{(m)}\Cal Sm_{(1)}
\otimes m_{(0)}$.
Note that $M$ is a left $H$-comodule via $\omega'$.
\end{example}

\begin{lemma}\label{star-omega.thm}
Let $\Cal M$ be a quasi-coherent $(G,\Cal O_X)$-module, and
$\frak n$ a quasi-coherent $\Cal O_X$-submodule of $\Cal M$.
Assume that the second projection $G\times X\rightarrow X$ is 
quasi-compact quasi-separated.
Then $\frak n^*$ agrees with the kernel of the composite map
\[
\Cal M\xrightarrow\omega a_*p_2^*\Cal M
\xrightarrow {\pi_{\frak n}} a_*p_2^*(\Cal M/\frak n),
\]
where $\pi_{\frak n}:\Cal M\rightarrow \Cal M/\frak n$ is the projection.
\end{lemma}

\proof The diagram
\[
\xymatrix{
 & \Cal M \ar[d]^u \ar[r]^-\omega
\ar `l [dddl] `[dddr]^{\id} [dddr] 
& a_*p_2^*\Cal M \ar[dd]^u \ar[r]^\pi 
& a_*p_2^*(\Cal M/\frak n) \ar[dd]^u \\
 & a_*a^*\Cal M \ar[ru]^\Phi \ar[d]^u & & \\
 & a_*E_*E^* a^*\Cal M\ar[r]^\Phi \ar[dr]^\cong &
a_*E_*E^* p_2^*\Cal M\ar[r]^-\pi \ar[d]^\cong  &
a_*E_*E^* p_2^*(\Cal M/\frak n) \ar[d]^\cong \\
 & & \Cal M \ar[r]^\pi & \Cal M/\frak n  \\
}
\]
is commutative, where $e:\Spec S\rightarrow G$ is the unit element, and
$E=e\times 1: X\rightarrow G\times X$.
Thus $\Ker(\pi\omega)\subset\Ker(u\pi\omega)=\Ker \pi=\frak n$.
Moreover, $\Ker(\pi\omega)\subset\Cal M$ is a quasi-coherent 
$(G,\Cal O_X)$-submodule of $\Cal M$, since $\pi\omega$ is 
$(G,\Cal O_X)$-linear.

So it suffices to show that 
any $(G,\Cal O_X)$-submodule $\Cal N$ of $
\Cal M$ contained in $\frak n$ is 
also
contained in $\Ker(\pi\omega)$.
This is trivial, since the diagram
\[
\xymatrix{
\Cal M \ar[r]^-\omega \ar[d]^\pi & a_*p_2^*\Cal M \ar[r]^-\pi \ar[d]^\pi & 
a_*p_2^*(\Cal M/\frak n) \\
\Cal M/\Cal N \ar[r]^-\omega & a_*p_2^*\Cal M/\Cal N \ar[ur]^-\pi
}
\]
is commutative.
\qed

\begin{definition}
We say that $X$ is $G$-integral (resp.\ $G$-reduced) 
if there is an integral (resp.\ reduced) closed subscheme
$Y$ of $X$ such that $Y^*=X$.
A $G$-ideal $\Cal P$ of $\Cal O_X$ is said to be $G$-prime (resp.\ 
$G$-radical), if
$V(\Cal P)$ is $G$-integral (resp.\ $G$-reduced).
\end{definition}

\begin{lemma}\label{Sim-G-int.thm}
Let $f:V\rightarrow X$ be a $G$-morphism of $G$-schemes.
If $V$ is $G$-integral \(resp.\ $G$-reduced\), 
then $\GSim(f)$ is $G$-integral \(resp.\ $G$-reduced\).
\end{lemma}

\proof There is an integral (resp.\ reduced) closed subscheme $Y$
of $V$ such that $Y^*=V$.
Then $Z:=\Sim(Y\hookrightarrow V\rightarrow X)$ is integral (resp.\ reduced),
see (\ref{reduced-sim.par}).
Then $\GSim(f)=Z^*$ by Lemma~\ref{GSim.thm}, and we are done.
\qed

\begin{corollary}
Let $f:V\rightarrow X$ be a $G$-morphism of $G$-schemes.
If $\I$ is a $G$-prime \(resp.\ $G$-radical\) ideal of $\O_V$,
then $(\I\cap \O_X)^*$ is $G$-prime \(resp.\ $G$-radical\).
\end{corollary}

\proof This is because $(\I\cap\O_X)^*$ defines $\GSim(V(\I)\hookrightarrow
V\rightarrow X)$, which is $G$-integral (resp.\ $G$-reduced).
\qed

\begin{lemma}\label{composition-GSim.thm}
Let
$f:W\rightarrow V$ and
$g:V\rightarrow X$ be $G$-morphisms of $G$-schemes, and
let $\iota:\GSim f\hookrightarrow V$ be the inclusion.
Then $\GSim(gf)=\GSim(g\iota)$.
\end{lemma}

\proof Similar to \cite[(9.5.5)]{EGA-I}.
\qed

\begin{lemma}
For a family $(\I_\lambda)$ of $G$-radical $G$-ideals of $\Cal O_X$,
$(\bigcap_\lambda\I_\lambda)^*$ is $G$-radical.
\end{lemma}

\proof There exists a family $(\J_\lambda)$ of radical ideals of $\Cal O_X$
such that $\J_\lambda^*=\I_\lambda$.
By 
Lemma~\ref{intersection.thm},
we have
\[
(\bigcap_\lambda\I_\lambda)^*=(\bigcap_\lambda\J_\lambda)^*
=((\bigcap_\lambda\J_\lambda)^\star)^*.
\]
By Lemma~\ref{radical-intersect.thm},
$(\bigcap_\lambda\J_\lambda)^\star$ is a radical ideal.
So $(\bigcap_\lambda\I_\lambda)^*$ is $G$-radical.
\qed

\begin{corollary}
The intersection of finitely many $G$-radical $G$-ideals is $G$-radical.
\qed
\end{corollary}

\begin{lemma}\label{principal-G-int.thm}
Let $Y$ be an $S$-scheme which is integral \(resp.\ reduced\).
Assume that the principal $G$-bundle $G\times Y$ is p-flat.
Then $G\times Y$ is $G$-integral \(resp.\ $G$-reduced\).
\end{lemma}

\proof Let us consider the closed subscheme $\{e\}\times Y$ of $G\times Y$,
where $e$ is the unit element.
Then
\[
G\times Y\cong G\times\{e\}\times Y\hookrightarrow G\times G\times Y
\xrightarrow{\mu\times 1}G\times Y
\]
is the identity, and its $G$-scheme theoretic image is $(\{e\}\times Y)^*$.
Hence $(\{e\}\times Y)^*=G\times Y$.
Since $\{e\}\times Y\cong Y$ is integral (resp.\ reduced), 
$G\times Y$ is $G$-integral (resp.\ $G$-reduced).
\qed

\begin{corollary}
Let $f:Y\rightarrow X$ be an $S$-morphism.
Assume that $Y$ is integral (resp.\ reduced).
Set $g$ to be the composite
\[
G\times Y \xrightarrow{1_G\times f}
G\times X \xrightarrow{a} X.
\]
Then $\GSim g$ is $G$-integral \(resp.\ $G$-reduced\).
\end{corollary}

\proof Follows from Lemma~\ref{principal-G-int.thm} and
Lemma~\ref{Sim-G-int.thm}.
\qed

\begin{lemma}\label{G-integral-open.thm}
Assume that the second projection $p_2:G\times X\rightarrow X$ is
quasi-compact quasi-separated.
If $X$ is $G$-integral \(resp.\ $G$-reduced\)
and $U$ is a 
non-empty
$G$-stable open subscheme of $X$, 
then $U$ is $G$-integral \(resp.\ $G$-reduced\).
\end{lemma}

\proof Take an integral closed subscheme $Y$ of $X$ such that $Y^*=X$.
Then
\[
\xymatrix{
G\times(Y\cap U) \ar[r]^-a \ar@{^{(}->}[d] & U \ar@{^{(}->}[d] \\
G\times Y \ar[r]^-a & X
}
\]
is a fiber square.
By (\ref{sim-flat.par}), $U=(Y\cap U)^*$.
Since $Y\cap U$ is an integral (reduced) 
closed subscheme of $U$, $U$ is $G$-integral (resp.\ $G$-reduced).
\qed

A morphism of schemes $f:Z\rightarrow Y$ is said to be {\em reduced} if
$f$ is flat, and 
for any point $y\in Y$ and any finite extension of fields $\kappa(y)
\rightarrow K$, the scheme $\Spec K\times_Y Z$ is reduced, where 
$\Spec K$ is viewed as a $Y$-scheme via the composite 
$\Spec K\rightarrow\Spec \kappa(y) \rightarrow Y$, see 
\cite[(6.8.1)]{EGA-IV}.
A base change of a reduced morphism is reduced.
The composite of two reduced morphisms is again reduced.
If $f:Z\rightarrow Y$ is a reduced morphism with $Y$ reduced,
then $Z$ is reduced.

\begin{lemma}
Assume that the second projection $p_2:G\times X\rightarrow X$ is
quasi-compact quasi-separated.
Let $\varphi:W\rightarrow X$ be a reduced $G$-morphism \(i.e., 
a $G$-morphism which is reduced as a morphism of schemes\) 
between $G$-schemes.
If $X$ is $G$-reduced, then $W$ is $G$-reduced.
\end{lemma}

\proof Similar to Lemma~\ref{G-integral-open.thm}.
\qed

\paragraph A $G$-ideal $\Cal P$ of $\Cal O_X$ is said to be 
{\em $G$-quasi-prime} if $\Cal P\neq\Cal O_X$, and if
$\I$ and $\J$ are $G$-ideals of $\Cal O_X$ such that $\I\J\subset\Cal P$,
then $\I\subset\Cal P$ or $\J\subset \P$ holds.

\begin{lemma}\label{G-prime-G-quasi-prime.thm}
If $\Cal P$ is a $G$-prime ideal of $\Cal O_X$, then
$\Cal P$ is $G$-quasi-prime.
\end{lemma}

\proof Let $\Cal P=\fp^*$ for a prime ideal $\fp$ of $\Cal O_X$.
Since $\Cal P\subset\fp\neq\Cal O_X$, $\Cal P\neq\Cal O_X$.
Let $\I$ and $\J$ be $G$-ideals of $\Cal O_X$ such that $\I\J\subset\P$.
Then $\I\J\subset\fp$.
Since $\fp$ is a prime ideal, we have $\I\subset\fp$ or $\J\subset\fp$.
Since $\I$ and $\J$ are $G$-ideals, we have
$\I\subset \fp^*=\Cal P$ or $\J\subset\fp^*=\Cal P$.
\qed

\paragraph
For a $G$-ideal $\Cal I\subset\Cal O_X$,
we denote the set of $G$-prime $G$-ideals containing $\Cal I$ by
$V_G(\Cal I)$.
The set $V_G(0)$ is denoted by $\Spec_G(X)$.
We define the {\em $G$-radical} of $\Cal I$ by
\[
\sqrt[G]{\Cal I}:=(\bigcap_{\Cal P\in V_G(\Cal I)}\Cal P)^*,
\]
where the right hand side is defined to be $\Cal O_X$ if 
$V_G(\Cal I)=\emptyset$.

\begin{definition}
For a $G$-ideal $\Cal I$, we define
\[
\Omega(\Cal I)=\{\Cal J\mid \text{$\Cal J$ is a quasi-cohernt $G$-ideal of
$\Cal O_X$, and } \Cal I\subset\Cal J\neq\Cal O_X\}.
\]
A maximal element of $\Omega(0)$ is 
said to be
{\em $G$-maximal}.
\end{definition}

\begin{lemma}\label{G-maximal-G-prime.thm}
A $G$-maximal $G$-ideal is $G$-quasi-prime.
If $X$ is quasi-compact, then a $G$-maximal $G$-ideal of $\Cal O_X$ 
is of the form $\fm^*$ for some maximal ideal $\fm$.
In particular, it is $G$-prime.
\end{lemma}

\proof We prove the first assertion.
Let $\Cal M$ be a $G$-maximal $G$-ideal of $\Cal O_X$,
and $\Cal I$ and $\Cal J$ be $G$-ideals such that
$\Cal I\Cal J\subset\Cal M$.
Assume that $\Cal I\not\subset \Cal M$ and $\Cal J\not\subset \Cal M$.
Then $\Cal I+\Cal M=\Cal O_X$ and $\Cal J+\Cal M=\Cal O_X$ by the 
$G$-maximality.
So
\[
\Cal O_X=(\Cal I+\Cal M)(\Cal J+\Cal M)\subset \Cal M+\Cal I\Cal J\subset
\Cal M.
\]
This contradicts $\Cal M\neq \Cal O_X$.

We prove the second assertion.
Let $\Cal M$ be a $G$-maximal $G$-ideal of $\Cal O_X$.
Since $X$ is quasi-compact, $V(\Cal M)$ is quasi-compact, and is 
clearly non-empty.
So there is a maximal ideal $\fm$ of $\Cal O_X$ containing $\Cal M$.
Since $\Cal M\subset\fm$, we have $\Cal M\subset \fm^*$.
By the maximality, $\Cal M=\fm^*$.
Since $\fm$ is a prime ideal, $\Cal M$ is 
$G$-prime.
\qed

\begin{lemma}\label{G-maximal.thm}
Let $\Cal I$ be a $G$-ideal of $\Cal O_X$.
If $\I\neq\Cal O_X$ and $X$ is quasi-compact,
then $\Omega(\Cal I)$ has a maximal element.
In particular, a non-empty quasi-compact p-flat $G$-scheme has a 
$G$-maximal $G$-ideal.
\end{lemma}

\proof Since $\Cal I\in\Omega(\Cal I)$, $\Omega(\Cal I)$ is non-empty.
So it suffices to show that, by Zorn's lemma, for any non-empty chain
(i.e., a totally ordered family) of elements $(\Cal J_\lambda)_{\lambda
\in\Lambda}$ of $\Omega(\Cal I)$, 
$\Cal J:=\sum_\lambda \Cal J_\lambda$ is again in $\Omega(\Cal I)$.
It is obvious that $\Cal J$ is a quasi-coherent $G$-ideal and
$\Cal J\supset \Cal I$.
It suffices to show that $\Cal J\neq \Cal O_X$.

Assume the contrary.
Let $X=\bigcup_{i=1}^n U_i$ be  a finite affine open covering of $X$,
which exists.
Let $J_{\lambda,i}:=\Gamma(U_i,\Cal J_\lambda)$, and $J_i=\Gamma(U_i,\Cal J)
=\sum_{\lambda}J_{\lambda,i}$.
As $1\in J_i$, there exists some $\mu_i$ such that $1\in J_{\mu_i,i}$.
When we set $\mu=\max(\mu_1,\ldots,\mu_n)$, then 
$\Cal J_\mu=\Cal O_X$.
This is a contradiction.
\qed

\begin{lemma}\label{trivial.thm}
Let $\Cal I$, $\Cal J$, $\Cal P$, and $\Cal I_\lambda$ 
$(\lambda\in\Lambda)$ be
$G$-ideals of $\Cal O_X$.
Then the following hold:
\begin{description}
\item[(i)] $\sqrt[G]{\Cal I}\supset \Cal I$.
\item[(ii)] If $\Cal I\supset \Cal J$, then $V_G(\Cal I)\subset V_G(\Cal J)$.
In particular, $\sqrt[G]{\Cal I}\supset \sqrt[G]{\Cal J}$.
\item[(iii)] $V_G(\sqrt[G]{\Cal I})=V_G(\Cal I)$.
In particular, $\sqrt[G]{\sqrt[G]{\Cal I}}=\sqrt[G]{\Cal I}$.
\item[(iv)] $V_G(\Cal I\Cal J)=V_G(\Cal I\cap \Cal J)=V_G(\Cal I)\cup
V_G(\Cal J)$.
So $\sqrt[G]{\Cal I\Cal J}=\sqrt[G]{\Cal I\cap\Cal J}=\grad{\Cal I}\cap
\grad{\Cal J}$.
\item[(v)] $V_G(\sum_\lambda \Cal I_\lambda)=
\bigcap_{\lambda}V_G(\Cal I_\lambda)$.
\item[(vi)] For $n\geq 1$, 
$\grad{\Cal I^n}=\grad{\Cal I}$.
\item[(vii)] If $\Cal P$ is $G$-prime, then $\grad{\Cal P}=\Cal P$.
\item[(viii)] If there exists some $n\geq 1$ such that $\Cal J^n
\subset \grad{\Cal I}$, then $\Cal J\subset \grad{\Cal I}$.
\end{description}
\end{lemma}

\proof The proof is easy, and left to the reader.
\qed

\begin{lemma} Let $\I$ be a $G$-ideal of $\Cal O_X$.
Then $\grad\I=\sqrt{\I}^{\,*}$.
In particular, 
$\I\subset\grad\I\subset\sqrt{\I}$ and hence
$\sqrt{\grad\I}=\sqrt\I$.
If, moreover, $X$ is noetherian, then 
there exists some $n\geq1$ such that $\grad\I^{\,n}\subset\Cal I$.
\end{lemma}

\proof In view of 
Lemma~\ref{trivial.thm}, {\bf (i)}, it suffices to prove the first assertion.
By Lemma~\ref{intersection.thm}, 
\[
\grad\I=(\bigcap_{\P\in V_G(\I)}\P)^*=
(\bigcap_{\fp\in V(\I)}\fp^*)^*
=
(\bigcap_{\fp\in V(\I)}\fp)^*=\sqrt{\I}^{\,*}.
\qed
\]

\begin{corollary}\label{G-radical-grad.thm}
Let $\I$ be a $G$-ideal of $\O_X$.
Then $\I$ is $G$-radical if and only if $\I=\grad\I$.
\end{corollary}

\proof Assume that $\I$ is $G$-radical so that $\I=\fa^*$ for a 
radical ideal $\fa$.
Since $\I=\fa^*\subset\fa$, $\sqrt\I\subset\sqrt\fa=\fa$.
So $\I\subset\grad\I=\sqrt{\I}^{\,*}\subset\fa^*=\I$, and
hence $\I=\grad\I$.

Conversely, assume that $\I=\grad\I$.
Then $\I=\grad\I=\sqrt{\I}^{\,*}$, and since $\sqrt{\I}$ is a radical ideal, 
$\I$ is $G$-radical.
\qed

\begin{corollary}\label{grad-O.thm}
For a $G$-ideal $\I$ of $\O_X$, $\I=\O_X$ if and only if $\grad{\I}=\O_X$
if and only if $V_G(\I)=\emptyset$.
\end{corollary}

\proof $V_G(\O_X)=\emptyset$ is trivial.
If $V_G(\I)=\emptyset$, then $\grad{\I}=\O_X$ by definition.
If $\grad{\I}=\O_X$, then $\sqrt{\I}=\O_X$, and hence $\I=\O_X$.
\qed

\begin{lemma}\label{sqrt.thm}
Assume that $X$ is noetherian.
If $\frak a$ is an ideal of $\Cal O_X$, then
$\grad{\frak a^*}=\sqrt{\frak a}^{\,*}$.
\end{lemma}

\proof Since $\grad{\frak a^*}\subset\sqrt{\frak a^*}\subset \sqrt{\frak a}$
and $\grad{\frak a^*}$ is a quasi-coherent $G$-ideal, 
we have $\grad{\frak a^*}\subset \sqrt{\frak a}^{\;*}$.

Next, set $\J:=\sqrt{\frak a}^{\;*}$.
Since $\J\subset\sqrt{\frak a}$, there exists some $n\geq1$
such that $\J^n\subset\frak a$.
Then $\J^n\subset \frak a^*\subset\grad{\frak a^*}$.
Hence $\J\subset \grad{\frak a^*}$.
\qed

\begin{definition}\label{G-point}
We say that $X$ is a {\em $G$-point} if $X$ has exactly two quasi-coherent
$G$-ideals of $\O_X$.
In other words, $X$ is a $G$-point if and only if $0$ is a $G$-maximal
$G$-ideal.
\end{definition}

\begin{corollary}
Let $X$ be a $G$-point, and $\M$ a quasi-coherent $(G,\O_X)$-module
of finite type.
Then $\M$ is locally free of a well-defined rank.
\end{corollary}

\proof 
Let $r$ be the smallest integer such that $\uFitt_r(\M)\neq 0$.
Since $\uFitt_r\M$ is a nonzero $G$-ideal, $\uFitt_r\M=\O_X$.
By Lemma~\ref{Fitting-locally-free.thm}, $\M$ is locally free of rank $r$.
\qed

\begin{example}
Let $k$ be a field, and $G$ a quasi-compact quasi-separated group scheme over $k$.
Then $X=G$ is a $G$-scheme with the left regular action.
Then the ascent functor $\pi^*:\Qch(\Spec k)\rightarrow
\Qch(G,X)$ is an equivalence by \cite[(12.10), (29.2)]{ETI}.
So it is easy to see that $X=G$ is a $G$-point.
In particular, $X$ is $G$-integral.
However, $G$ may not be reduced or irreducible, even if $G$ is 
affine and finite type.
So a $G$-integral $G$-scheme may not be reduced or irreducible in general.
\end{example}

\section{$G$-primary $G$-ideals}

\paragraph As in the last section, let $S$ be a scheme, $G$ an $S$-group
scheme, and $X$ a p-flat $G$-scheme.
In this section, we always assume that $X$ is noetherian.
Let $\M$ be a coherent $(G,\O_X)$-module.
Unless otherwise specified, a submodule of a coherent $\Cal O_X$-module
means a coherent submodule.
In particular, an ideal of $\O_X$ means a coherent ideal.
Let $\N$ be a (coherent) $(G,\O_X)$-submodule of $\M$.

\begin{lemma}\label{prime-prelim.thm}
Let $\P$ be a $G$-quasi-prime $G$-ideal of $\O_X$.
Then $\P=\sqrt{\P}^{\,*}$.
\end{lemma}

\proof Since $\P\subset\sqrt{\P}$, we have
$\P\subset\sqrt{\P}^{\,*}$.

Now we set $\J:=\sqrt{\P}^{\,*}$ and we prove $\J\subset\P$.
Since $X$ is noetherian, there exists some $n\geq 1$ such that
$\sqrt{\P}^{\,n}\subset\P$.
Then $\J^n\subset\sqrt{\P}^{\,n}\subset\P$.
Since $\P$ is $G$-quasi-prime, $\J\subset\P$, as desired.
\qed

\begin{lemma}\label{G-prime.thm}
Let $\P$ be a $G$-ideal of $\O_X$.
Then $\P$ is a $G$-prime $G$-ideal if and only if $\P$ is $G$-quasi-prime.
If this is the case, $\P=\fp^*$ for some minimal prime ideal $\fp$ of $\P$.
\end{lemma}

\proof The `only if' part is Lemma~\ref{G-prime-G-quasi-prime.thm}.

Assume that $\P$ is $G$-quasi-prime.
Let
\[
\sqrt{\Cal P}=\frak p_1\cap\cdots\cap \frak p_r
\]
be a minimal prime decomposition so that each $\frak p_i$ is a minimal 
prime ideal of $\Cal P$.
Then by Lemma~\ref{prime-prelim.thm},
\[
\Cal P=\sqrt{\Cal P}^{\,*}
=\frak p_1^*\cap\cdots\cap\frak p_r^*.
\]
As $\Cal P$ is $G$-quasi-prime, there exists some $i$ such that
$\frak p_i^*\subset\Cal P$.
So $\frak p_i^*=\Cal P$.
\qed

\begin{definition}
A $(G,\cO_X)$-submodule $\Cal N$ is said to be a {\em $G$-primary}
submodule of $\Cal M$ if
\begin{description}
\item[(i)] $\Cal M\neq\Cal N$, and
\item[(ii)] For a $(G,\Cal O_X)$-submodule $\Cal L$ of $\Cal M$ and
for a $G$-ideal $\I$ of $\cO_X$, if 
$\I\Cal L\subset \Cal N$ and $\Cal L\not\subset\Cal N$, then
$\I\subset\grad{\Cal N:\Cal M}$.
\end{description}
A $G$-primary $(G,\O_X)$-submodule 
of $\O_X$ is said to be a $G$-primary $G$-ideal.
If $0$ is a $G$-primary submodule of $\Cal M$, then we say that $\Cal M$
is a $G$-primary $(G,\Cal O_X)$-module.
If $\Cal O_X$ is a $G$-primary module, 
then we say that $X$ is $G$-primary.
\end{definition}

\begin{lemma}\label{G-primary-G-prime.thm}
Let $\Cal Q$ be a $G$-primary ideal of $\Cal O_X$.
Then $\grad\Q$ is $G$-prime.
\end{lemma}

\proof It suffices to prove that $\grad\Q$ is $G$-quasi-prime.

Since $\Q\neq\cO_X$, we have $\grad\Q\neq\cO_X$ by
Corollary~\ref{grad-O.thm}.

Let $\I$ and $\J$ be quasi-coherent $G$-ideals.
Assume that $\Cal I\Cal J\subset\grad\Q$ and $\J\not\subset\grad\Q$.
Then there exists some $n\geq 1$ such that 
$\Cal I^n\Cal J^n\subset \Cal Q$ and $\Cal J^n\not\subset\grad \Q$.
Hence $\I^n\subset\Q$.
This shows that $\I\subset\grad\Q$.
\qed

\begin{lemma}\label{module-ideal-G-primary.thm}
Let $\Cal N$ be a $G$-primary coherent $(G,\cO_X)$-submodule
of $\Cal M$.
Then $\Cal N:\Cal M$ is $G$-primary.
In particular, $\grad{\Cal N:\Cal M}$ is $G$-prime.
\end{lemma}

\proof Let $\I$ and $\J$ be coherent $G$-ideals, and assume that
$\I\J\subset\Cal N:\Cal M$ and that $\J\not\subset\Cal N:\Cal M$.
Then $\J\Cal M\not\subset\Cal N$ and $\I\J\Cal M\subset\Cal N$.
Hence $\I\subset \grad{\Cal N:\Cal M}$.
The last assertion follows from
Lemma~\ref{G-primary-G-prime.thm}.
\qed

\paragraph If $\Cal N$ is $G$-primary and $\grad{\Cal N:\Cal M}=\Cal P$,
then we say that $\Cal N$ is $\Cal P$-$G$-primary.

\begin{lemma}\label{primary.thm}
If $\frak m$ is a $\frak p$-primary submodule of $\Cal M$, then
$\frak m^*$ is a $\frak p^*$-$G$-primary $(G,\Cal O_X)$-submodule of $\M$.
\end{lemma}

\proof By Lemma~\ref{sqrt.thm} and Lemma~\ref{colon.thm},
\[
\frak p^*=
\sqrt{\frak m:\Cal M}^{\;*}
=\grad{(\frak m:\Cal M)^{\;*}}
=\grad{\frak m^*:\Cal M}.
\]

Let $\Cal I$ be a coherent $G$-ideal of $\cO_X$ and
$\Cal L$ a coherent $(G,\Cal O_X)$-submodule of $\Cal M$.
Assume that $\Cal I\Cal L\subset \frak m^*$ and $\Cal L\not\subset\frak m^*$.
Then $\Cal I\Cal L\subset \frak m$ and $\Cal L\not\subset\frak m$.
Hence $\Cal I\subset\sqrt{\frak m:\Cal M}$.
Hence 
\[
\Cal I\subset
\sqrt{\frak m:\Cal M}^{\;*}=\frak p^*=
\grad{\frak m^*:\Cal M}.
\]
Since $\frak m\neq\Cal M$, we have $\frak m^*\neq\Cal M$.
\qed

\begin{lemma}\label{joining.thm}
Let $\Cal N$ and $\Cal N'$ be $\Cal P$-$G$-primary coherent 
$(G,\Cal O_X)$-submodules of $\Cal M$.
Then $\Cal N\cap \Cal N'$ is also $\Cal P$-$G$-primary.
\end{lemma}

\proof $(\Cal N:\Cal M)\cap(\Cal N':\Cal M)=(\Cal N\cap \Cal N'):\Cal M$.
So
\begin{multline*}
\grad{(\Cal N\cap \Cal N'):\Cal M}=
\grad{(\Cal N:\Cal M)\cap(\Cal N':\Cal M)}\\
=\grad{\Cal N:\Cal M}\cap\grad{\Cal N':\Cal M}
=\Cal P\cap\Cal P=\Cal P.
\end{multline*}
Let $\Cal I$ be a coherent $G$-ideal, and $\Cal L$ a coherent 
$(G,\Cal O_X)$-submodule of $\Cal M$ such that
$\Cal I\Cal L\subset\Cal N\cap\Cal N'$ and $\Cal L\not\subset\Cal N
\cap\Cal N'$.
Then $\Cal L\not\subset\Cal N$ or $\Cal L\not\subset\Cal N'$.
If $\Cal L\not\subset\Cal N$, then $\Cal I\subset
\grad{\Cal N:\Cal M}=\Cal P$.
If $\Cal L\not\subset\Cal N'$, then $\Cal I\subset
\grad{\Cal N':\Cal M}=\Cal P$.
\qed

\begin{definition}
Let $\Cal M$ be a coherent $(G,\Cal O_X)$-module, and
$\Cal N$ a coherent $(G,\Cal O_X)$-submodule of $\Cal M$.
An expression
\begin{equation}\label{G-primary.eq}
\Cal N=\Cal N_1\cap\cdots\cap \Cal N_r
\end{equation}
is called a {\em $G$-primary decomposition} if this equation holds, and
each $\Cal N_i$ is $G$-primary.
The $G$-primary decomposition (\ref{G-primary.eq}) is said to be {\em
irredundant} if for each $i$, $\bigcap_{j\neq i}\Cal N_j\neq \Cal N$.
It is said to be {\em minimal} if it is irredundant and 
$\grad{\Cal N_i:\Cal M}\neq \grad{\Cal N_j:\Cal M}$ for $i\neq j$.
\end{definition}

\begin{lemma}\label{existence.thm}
Let $\Cal M$ be a coherent $(G,\Cal O_X)$-module, and
$\Cal N$ a coherent $(G,\Cal O_X)$-submodule of $\Cal M$.
Then $\Cal N$ has a minimal $G$-primary decomposition.
\end{lemma}

\proof Let $\Cal N=\frak n_1\cap \cdots\cap\frak n_s$ be a primary 
decomposition, which exists by (\ref{primary.par}).
Then
\[
\Cal N=\Cal N^*=(\frak n_1\cap \cdots\cap\frak n_s)^*
=\frak n_1^*\cap \cdots\cap\frak n_s^*.
\]
This is a $G$-primary decomposition by Lemma~\ref{primary.thm}.
Omitting redundant terms, we get an irredundant decomposition.
Let 
\[
\Cal N=\Cal N_1\cap\cdots\cap \Cal N_t
\]
be the decomposition so obtained.
We say that $i\sim j$ if $\grad{\Cal N_i:\Cal M}=\grad{\Cal N_j:\Cal M}$.
Let $E_1,\ldots,E_r$ be the equivalence classes with respect to the
equivalence relation $\sim$.
Then letting $\Cal E_i=\bigcap_{j\in E_i}\Cal N_j$,
\[
\Cal N=\Cal E_1\cap\cdots\cap \Cal E_r
\]
is a minimal $G$-primary decomposition by 
Lemma~\ref{joining.thm}.
\qed

\begin{lemma}\label{key-lemma.thm}
Let $\Cal M$ be a coherent $(G,\Cal O_X)$-module,
$\Cal N$ a coherent $(G,\Cal O_X)$-submodule of $\Cal M$, and
$\Cal P$ a coherent $G$-ideal of $\Cal O_X$.
Then the following are equivalent.
\begin{description}
\item[(i)] $\Cal N$ is $\Cal P$-$G$-primary.
\item[(ii)] The following three conditions hold:
\begin{description}
\item[(a)] $\Cal N\neq \Cal M$.
\item[(b)] $\Cal P\subset\grad{\Cal N:\Cal M}$.
\item[(c)] If $\Cal L$ is a coherent $(G,\Cal O_X)$-submodule of $\Cal M$,
$\Cal J$ is a $G$-ideal, $\Cal L\not\subset \Cal N$, 
and $\J\not\subset\Cal P$, then $\J\Cal L\not\subset\Cal N$.
\end{description}
\end{description}
\end{lemma}

\proof {\bf (i)$\Rightarrow$(ii)} is clear.

{\bf (ii)$\Rightarrow$(i)} Set $\Cal K:=\grad{\Cal N:\Cal M}$.
We show $\Cal P=\Cal K$.
There exists some $n\geq 1$ $\Cal K^n\Cal M\subset\Cal N$ and
$\Cal K^{n-1}\Cal M\not\subset\Cal N$ by {\bf (a)}.
Then by {\bf (c)}, $\Cal K\subset\Cal P$.
So $\Cal K=\Cal P$ by {\bf (b)}.
By {\bf (a)} and {\bf (c)}, we have that $\Cal N$ is $\Cal P$-$G$-primary.
\qed

\begin{lemma}\label{colon2.thm}
Let $\Cal M$ be a coherent $(G,\Cal O_X)$-module, and
$\Cal N$ a $G$-primary coherent $(G,\Cal O_X)$-submodule of $\Cal M$.
Set $\Cal P=\grad{\Cal N:\Cal M}$.
Let $\Cal L$ be a coherent $(G,\Cal O_X)$-submodule of $\Cal M$,
and $\Cal I$ a coherent $G$-ideal of $\Cal O_X$.
Then the following hold.
\begin{description}
\item[(i)] If $\Cal L\subset\Cal N$, then $\Cal N:\Cal L=\Cal O_X$.
\item[(ii)] If $\Cal L\not\subset\Cal N$, then $\Cal N:\Cal L$ is
$\Cal P$-$G$-primary.
\item[(iii)] If $\Cal I\subset \Cal N:\Cal M$, then $\Cal N:\Cal I=\Cal M$.
\item[(iv)] If $\Cal I\not\subset\Cal N:\Cal M$, then 
$\Cal N:\Cal I$ is $\Cal P$-$G$-primary.
\item[(v)] If $\Cal I\not\subset\Cal P$, then $\Cal N:\I=\Cal N$.
\end{description}
\end{lemma}

\proof {\bf (i)} and {\bf (iii)} are trivial.

{\bf (ii)} $\Cal N:\Cal L\neq\Cal O_X$ is trivial.
$\grad{\Cal N:\Cal L}\supset\grad{\Cal N:\Cal M}=\Cal P$.
Let $\Cal J$ and $\Cal K$ be coherent $G$-ideals of $\Cal O_X$ such that
$\Cal J\not\subset\Cal N:\Cal L$ and $\Cal K\not\subset\Cal P$.
Then $\Cal J\Cal L\not\subset\Cal N$, and
$\Cal K\not\subset\Cal P$.
So $\Cal J\Cal K\Cal L\not\subset\Cal N$.
This shows $\Cal J\Cal K\not\subset\Cal N:\Cal L$.
By Lemma~\ref{key-lemma.thm}, $\Cal N:\Cal L$ is $\Cal P$-$G$-primary.

{\bf (iv)} Since $\I\Cal M\not\subset\Cal N$, $\Cal M\neq\Cal N:\I$.
We have
\[
\grad{(\Cal N:\Cal I):\Cal M}
=
\grad{\Cal N:\Cal I\Cal M}
\supset
\grad{\Cal N:\Cal M}
=
\Cal P.
\]
Let $\Cal L$ be a coherent $(G,\Cal O_X)$-submodule of $\Cal M$, and
$\Cal J$ be a coherent $G$-ideal of $\Cal O_X$ such that
$\Cal L\not\subset\Cal N:\Cal I$ and $\Cal J\not\subset\Cal P$.
Since $\Cal I\Cal L\not\subset\Cal N$ and 
$\Cal J\not\subset\Cal P$, we have that
$\Cal I\Cal J\Cal L\not\subset\Cal N$.
This shows $\Cal J\Cal L\not\subset\Cal N:\Cal I$.
By Lemma~\ref{key-lemma.thm}, 
$\Cal N:\Cal I$ is $\Cal P$-$G$-primary.

{\bf (v)} If $\Cal N:\Cal I\not\subset\Cal N$, then as $\Cal I\not\subset
\Cal P$, we have that $\I(\Cal N:\I)\not\subset\Cal N$.
This is a contradiction.
\qed

\paragraph Let $\Cal M$ be a coherent $(G,\Cal O_X)$-module, 
and $\Cal N$ a coherent $(G,\Cal O_X)$-submodule of $\Cal M$.
Let 
\begin{equation}\label{decomposition.eq}
\Cal N=\Cal Q_1\cap\cdots\cap\Cal Q_r
\end{equation}
be a minimal $G$-primary decomposition, which exists by 
Lemma~\ref{existence.thm}.
Set $\Cal M_i=\bigcap_{j\neq i}\Cal Q_j$, and
$\Cal P_i=\grad{\Cal Q_i:\Cal M}$.

\begin{theorem}\label{G-ass.thm}
We have
\begin{multline*}
\{\Cal P_1,\ldots,\Cal P_r\}=
\{\Cal N:\Cal L\mid
\text{$\Cal L$ is a coherent $(G,\Cal O_X)$-submodule of $\Cal M$, }
\\
\text{and $\Cal N:\Cal L$ is $G$-prime}\}.
\end{multline*}
In particular, this set depends only on $\Cal M/\Cal N$, and independent of the
choice of minimal $G$-primary decomposition of $\Cal N$.
\end{theorem}

\proof Since the decomposition (\ref{decomposition.eq}) is irredundant,
\[
\Cal N:\Cal M_i=\bigcap_{j=1}^r(\Cal Q_j:\Cal M_i)=
\Cal Q_i:\Cal M_i\neq\Cal O_X.
\]
Thus $\Cal N:\Cal M_i$ is $\Cal P_i$-$G$-primary by Lemma~\ref{colon2.thm},
{\bf (ii)}.
Take the minimum $n\geq 1$ such that $\Cal P_i^n\subset \Cal N:\Cal M_i$,
and set $\Cal L:=\Cal P_i^{n-1}\Cal M_i$.
Since $\Cal P_i^{n-1}\not\subset \Cal N:\Cal M_i$, 
$\Cal N:\Cal L=\Cal N:\Cal P_i^{n-1}\Cal M_i=(\Cal N:\Cal M_i):\Cal P_i^{n-1}
$ is $\Cal P_i$-$G$-primary.
In particular, 
$\Cal N:\Cal L\subset\grad{\Cal N:\Cal L}=\Cal P_i$.
On the other hand,
\[
(\Cal N:\Cal L):\Cal P_i
=
\Cal N:\Cal P^n_i\Cal M_i
=
(\Cal N:\Cal M_i):\Cal P_i^n=\Cal O_X.
\]
Thus $\Cal P_i\subset \Cal N:\Cal L$.
Hence $\Cal N:\Cal L=\Cal P_i$.
Thus each $\Cal P_i$ is of the form $\Cal N:\Cal L$ for some $\Cal L$.

Conversely, let $\Cal L\subset\Cal M$, and assume that
$\Cal N:\Cal L$ is $G$-prime.
Set $\Cal P=\Cal N:\Cal L$.
We show that $\Cal P=\Cal P_i$ for some $i$.
Renumbering if necessary, we may assume that $\Cal L\not\subset
\Cal Q_i$ if and only if $i\leq s$.
Then by Lemma~\ref{colon2.thm}, 
{\bf (ii)}, 
\[
\Cal P=\grad{\Cal P}=\grad{\Cal N:\Cal L}
=
\grad{\bigcap_{i=1}^s(\Cal Q_i:\Cal L)}
=
\bigcap_{i=1}^s\grad{\Cal Q_i:\Cal L}
=
\bigcap_{i=1}^s\Cal P_i.
\]
So $s\geq 1$, and there exists some $i$ such that $\Cal P_i=\Cal P$.
\qed

\begin{definition}
We set $\Ass_G(\Cal M/\Cal N)=\{\Cal P_1,\ldots,\Cal P_r\}$.
Note that $\Ass_G(\Cal M/\Cal N)$ depends only on $\Cal M/\Cal N$.
An element of $\Ass_G(\Cal M/\Cal N)$ is called a $G$-associated $G$-prime
$G$-ideal
of $\Cal M/\Cal N$ (however, also called a $G$-associated $G$-prime of
the submodule $\Cal N$).
The set of minimal elements in $\Ass_G(\Cal M/\Cal N)$ is denoted by
$\Min_G(\Cal M/\Cal N)$.
An element of $\Min_G(\Cal M/\Cal N)$ is called a minimal $G$-prime 
$G$-ideal of $\Cal M/\Cal N$.
An element of $\Ass_G(\Cal M/\Cal N)\setminus \Min_G(\Cal M/\Cal N)$ is 
called an embedded $G$-prime $G$-ideal.
\end{definition}

\begin{proposition}\label{uniqueness.thm}
Let $\Omega$ be a subset of $\{\Cal P_1,\ldots,\Cal P_r\}$.
Assume that $\Cal P_i\subset \Cal P_j\in \Omega$ implies $\Cal P_i
\in \Omega$ for $1\leq i,j\leq r$.
Then $\bigcap_{\Cal P_j\in\Omega}\Cal Q_j$ is independent of the choice of
minimal $G$-primary decomposition.
\end{proposition}

\proof Set $\J:=\bigcap_{\Cal P_i\notin\Omega}\Cal P_i$.
It suffices to prove that $\bigcap_{\Cal P_j\in\Omega}\Cal Q_j=
\Cal N:\J^n$ for $n\gg 0$.

For $\Cal P_i\notin\Omega$, $\J\subset\Cal P_i=\grad{\Cal Q_i:\Cal M}$.
Hence there exists some $n_0$ such that $\J^{n_0}\subset\Cal Q_i:\Cal M$ for
all $i$ such that $\Cal P_i\notin\Omega$.
Take $n$ so that $n\geq n_0$.
Then $\Cal Q_i:\J^n=\Cal M$, since $\J^n\Cal M\subset\Cal Q_i$, for 
$\Cal P_i\notin\Omega$.

If $\Cal P_i\notin\Omega$ and $\Cal P_j\in\Omega$, then
$\Cal P_i\not\subset\Cal P_j$ by assumption.
Hence $\J\not\subset\Cal P_j$.
Hence $\J^n\not\subset\Cal P_j$.
Thus $\Cal Q_j:\J^n=\Cal Q_j$ for $\Cal P_j\in\Omega$ by
Lemma~\ref{colon2.thm}, {\bf(v)}.

So $\Cal N:\J^n=\bigcap_i\Cal Q_i:\J^n=\bigcap_{\Cal P_j\in\Omega}\Cal Q_j$.
\qed

\begin{lemma}\label{Ass-Ass.thm}
Let $\Cal M$ be a coherent $(G,\Cal O_X)$-module, and
$\Cal N$ a coherent $(G,\Cal O_X)$-submodule of $\Cal M$.
If $\Cal P\in\Ass_G(\Cal M/\Cal N)$ and
$\frak p\in\Ass(\Cal O_X/\Cal P)$, then
$\frak p\in\Ass(\Cal M/\Cal N)$.
\end{lemma}

\proof There exists some coherent $(G,\Cal O_X)$-submodule $\Cal L$ of 
$\Cal M$ such that $\Cal N:\Cal L=\Cal P$.
Moreover, there exists some coherent ideal $\frak a$ of $\Cal O_X$ such
that $\frak p=\Cal P:\frak a$.
Then
\[
\frak p=\Cal P:\frak a=(\Cal N:\Cal L):\frak a=\Cal N:\frak a\Cal L.
\]
Thus $\frak p\in\Ass(\Cal M/\Cal N)$.
\qed

\begin{lemma}\label{min-prime.thm}
Let $\Cal M$ be a coherent $(G,\Cal O_X)$-module, and
$\Cal N$ a coherent $(G,\Cal O_X)$-submodule of $\Cal M$.
Then we have $V_G(\N:\M)=\bigcup_{\P\in\Min_G(\M/\N)}V_G(\P)$.
In particular, for a quasi-coherent $G$-prime 
$G$-ideal 
$\Cal P$ of $\Cal O_X$, 
$\Cal P\in\Min_G(\Cal M/\Cal N)$ if and only if $\Cal P$ is a minimal
element of $V_G(\Cal N:\Cal M)$.
Moreover, we have 
$\bigcap_{\Q\in V_G(\N:\M)}\Q=\bigcap_{\P\in\Min_G(\M/\N)}\P$.
\end{lemma}

\proof We may assume that $\Cal N=0$.
Let 
\[
0=\Cal Q_1\cap\cdots\cap\Cal Q_r
\]
be a minimal $G$-primary decomposition of $0$ in $\Cal M$.
Set $\Cal P_i:=\grad{\Cal Q_i:\Cal M}$.
Then
\[
\Cal P_i\supset \Cal Q_i:\Cal M\supset 0:\Cal M.
\]
In particular, $V_G(\Cal P_i)\subset V_G(0:\Cal M)$.
On the other hand, $\Cal M$ is a submodule of $\bigoplus_{i=1}^r\Cal M/\Cal 
Q_i$.
So $0:\Cal M\supset 0:\bigoplus_i \Cal M/\Cal Q_i$.
So
\[
V_G(0:\Cal M)\subset V_G(0:\bigoplus_i\Cal M/\Cal Q_i)
=V_G(\bigcap_i(\Cal Q_i:\Cal M))=\bigcup_i V_G(\Cal P_i)\subset V_G(0:\Cal M).
\]
So $V_G(0:\Cal M)=\bigcup_i V_G(\Cal P_i)$.
If $\P_i$ is an embedded $G$-prime $G$-ideal, then $V_G(\P_i)$ in the union is
redundant, and we have $V_G(0:\M)=\bigcup_{\P\in\Min_G \M}V_G(\P)$.
The rest of the assertions are now obvious.
\qed

\begin{corollary}\label{star-unnecessary.thm}
For an ideal $\I$ of $\Cal O_X$, $\sqrt[G]{\I}=\bigcap_{\cP\in V_G(\I)}\cP
=\bigcap_{\cP\in\Min_G(\cO_X/\I)}\cP$.
\end{corollary}

\proof 
By Lemma~\ref{min-prime.thm}, 
\[
\sqrt[G]{\I}=(\bigcap_{\cP\in V_G(\I)}\cP)^*=
(\bigcap_{\cP\in\Min_G(\cO_X/\I)}\cP)^*
=\bigcap_{\cP\in\Min_G(\cO_X/\I)}\cP
=\bigcap_{\cP\in V_G(\I)}\cP,
\]
since $\Min_G(\cO_X/\I)$ is a finite set and the intersection of finitely
many quasi-coherent $G$-ideals is a quasi-coherent $G$-ideal.
\qed

\begin{lemma}\label{primary-star.thm}
Let $\Cal M$ be a coherent $(G,\Cal O_X)$-module, and $\Cal N$ a 
$G$-primary coherent $(G,\Cal O_X)$-submodule.
Let 
\[
\Cal N=\frak n_1\cap\cdots\cap \frak n_r
\]
be a minimal primary decomposition.
If $\sqrt{\frak n_1:\Cal M}$ is a minimal prime ideal, 
then $\frak n_1^*=\Cal N$.
In particular, $\sqrt{\fn_1:\M}^{\,*}=\grad{\N:\M}$.
\end{lemma}

\proof 
First assume that $\N:\fn_1^*\subset\grad{\N:\M}$.
Then
\[
\bigcap_{i\geq 2}(\fn_i:\M)\subset\bigcap_{i\geq 1}(\fn_i:\fn_1^*)
=\N:\fn_1^*\subset \grad{\N:\M}\subset\sqrt{\N:\M}\subset\sqrt{\fn_1:\M}.
\]
This contradicts the minimality of $\sqrt{\fn_1:\M}$.
Hence $\N:\fn_1^*\not\subset\grad{\N:\M}$.
Since $\fn_1^*(\N:\fn_1^*)\subset \Cal N$ and $\Cal N$ is
$G$-primary, $\frak n_1^*\subset \Cal N$.
As $\frak n_1^*\supset \Cal N$ is trivial, $\frak n_1^*=\Cal N$.

Hence,
\[
\sqrt{\fn_1:\M}^{\,*}=\grad{(\fn_1:\M)^*}=\grad{\fn_1^*:\M}=\grad{\N:\M}.
\]
\qed

\begin{corollary}\label{G-prime-minimal}
Let $\P$ be a $G$-prime $G$-ideal.
For any minimal prime ideal $\fp$ of $\P$, we have $\fp^*=\P$.
\end{corollary}

\proof Let $\P=\fq_1\cap\cdots\cap\fq_r$ be a minimal primary decomposition
such that $\sqrt{\fq_1}=\fp$.
Then, $\fp^*=\sqrt{\fq_1}^{\,*}=\grad{\fq_1^*}=\grad{\P}=\P$.
\qed

\begin{corollary}
Let $\M$ be a coherent $(G,\O_X)$-module, 
and $\N$ a $(G,\O_X)$-submodule.
Then $\N$ is a $G$-primary submodule of $\M$
if and only if $\N=\fn^*$ for some
primary submodule $\fn$ of 
$\M$.
\end{corollary}

\proof The `if' part is Lemma~\ref{primary.thm}.
The `only if' part follows from Lemma~\ref{primary-star.thm}.
\qed

\begin{lemma}
For a $G$-ideal $\Cal I$ of $\Cal O_X$,
the following are equivalent.
\begin{description}
\item[(i)] $\Cal I$ is $G$-radical.
\item[(ii)] There are finitely many $G$-prime $G$-ideals
$\Cal P_1,\ldots,\Cal P_r$ of $\Cal O_X$ such that
$\Cal I=\Cal P_1\cap\cdots\cap \Cal P_r$.
\end{description}
\end{lemma}

\proof {\bf (i)$\Rightarrow$(ii)} 
follows from Corollary~\ref{star-unnecessary.thm}
and
Corollary~\ref{G-radical-grad.thm}.

{\bf (ii)$\Rightarrow$(i)} 
$\sqrt[G]{\Cal I}=\sqrt[G]{\Cal P_1\cap\cdots\cap\Cal P_r}
=\sqrt[G]{\Cal P_1}\cap\cdots\cap\sqrt[G]{\Cal P_r}
=\Cal P_1\cap\cdots\cap\Cal P_r=\Cal I$.
\qed

\section{Group schemes of finite type}

In this section, $S$ is a scheme, $G$ an $S$-group scheme,
$X$ a p-flat noetherian $G$-scheme, and $\Cal M$ a coherent
$(G,\Cal O_X)$-module.
In this section, we assume that $p_2:G\times X\rightarrow X$ is
of finite type.

Let $\Cal N$ be a coherent $(G,\Cal O_X)$-submodule of $\Cal M$.

\begin{lemma} Let
\[
\N=\fn_1\cap\cdots\cap\fn_r\cap\fn_{r+1}\cap\cdots\cap \fn_{r+s}
\]
be a minimal primary decomposition such that
\[
\Min(\M/\N)=\{\fp_1,\ldots,\fp_r\},
\]
where $\fp_i=
\sqrt{\fn_i:\M}$.
Then
\begin{description}
\item[(i)] The $(S_1)$-locus of $\Cal M/\Cal N$ is
$X\setminus\bigcup_{i=1}^s\Supp(\Cal M/\fn_{i+r})$.
\item[(ii)] The $(S_1)$-locus of $\M/\N$ is a $G$-stable open subset of $X$.
\item[(iii)] $\fn_1\cap\cdots\cap 
\fn_r$ is a coherent $(G,\Cal O_X)$-submodule of
$\M$.
\end{description}
\end{lemma}

\proof Replacing $\M$ by $\M/\N$, we may assume that $\N=0$.
Since $X\setminus\Supp\M$ is $G$-stable open, replacing $X$ by $V(\uann \M)$, 
we may assume that $\uann\M=0$.

{\bf (i)} 
Note that a coherent $\Cal O_X$-module $\Cal L$ satisfies Serre's
$(S_1)$-condition at $x\in X$ if and only if $\Cal L_x$ does not have
an embedded prime ideal.
The assertion follows from this.

{\bf (ii)} 
Let $U$ be the $(S_1)$ locus 
$X\setminus\bigcup_{i=1}^s\Supp(\M/\fn_{i+r})$ of $\M$.
It is an open subset.
Since the action $a:G\times X\rightarrow X$ and the second projection
$p_2:G\times X\rightarrow X$ are Cohen--Macaulay morphisms by
\cite[Lemma~31.14]{ETI}, both $a^{-1}(U)$ and $p_2^{-1}(U)=G\times U$
are the $(S_1)$-locus of $a^*\Cal M\cong p_2^*\Cal M$ by 
\cite[(6.4.1)]{EGA-IV}.
So $a^{-1}(U)=G\times U$, and $U$ is $G$-stable.

{\bf (iii)} 
Let $\iota:U\hookrightarrow X$ be the inclusion.
It suffices to show that $\uGamma_{X,U}(\Cal M):=\Ker(\Cal M
\rightarrow \iota_*\iota^*\Cal M)$ agrees with $\fn_1\cap\cdots\cap\fn_r$,
see for the notation, \cite[(3.1)]{LCDS}.

Since the composite
\[
\fn_1\cap\cdots\cap\fn_r\hookrightarrow\M
\rightarrow\bigoplus_{i=1}^s\M/\fn_{i+r}
\]
is 
a 
mono, $\iota^*(\fn_1\cap\cdots\cap\fn_r)=0$.
Hence $\fn_1\cap\cdots\cap\fn_r\subset\uGamma_{X,U}(\Cal M)$.
It suffices to show that $\uGamma_{X,U}(\M/\fn_1\cap\cdots\cap\fn_r)=0$.
(If so, $\fn_1\cap\cdots\cap\fn_r=\uGamma_{X,U}(
\fn_1\cap\cdots\cap\fn_r)=
\uGamma_{X,U}(\Cal M)$).
As $
\M/\fn_1\cap\cdots\cap\fn_r\subset\bigoplus_{i=1}^r\M/\fn_i$,
it suffices to show that $\uGamma_{X,U}(\M/\fn_i)=0$ for $i\leq r$.
Assume the contrary.
Then since $\Ass(\uGamma_{X,U}(\M/\fn_i))\subset\Ass(\M/\fn_i)$,
$\Ass(\M/\fn_i)$ contains a point in $X\setminus U$.
On the other hand, $\Ass(\M/\fn_i)$ is a singleton, and its point is
a generic point of an irreducible component of $X$.
As $U$ is dense in $X$, this is a contradiction.
\qed

\begin{corollary}\label{G-primary-no-emb.thm}
If $\Cal N$ is $G$-primary, then $\M/\N$ does not have an embedded prime
ideal.
\end{corollary}

\proof We may assume that $\N=0$ and $\uann\M=0$.
Let
\[
0=\fn_1\cap\cdots\cap\fn_r\cap\fn_{r+1}\cap\cdots\cap\fn_{r+s}
\]
be a minimal primary decomposition such that
$\Min(\Cal M)=\{\fp_1,\ldots,\fp_r\}$, where $\fp_i=\sqrt{\fn_i:\Cal M}$.
Set $\Cal L=\fn_1\cap\cdots\cap\fn_r$.
It suffices to show that $\Cal L=0$.
Set $\J:=\uann\Cal L$.
Note that $\Cal L$ is a coherent $(G,\Cal O_X)$-submodule of $\Cal M$ by
the lemma, and $\J$ is a coherent $G$-ideal.
Since $V(\J)\subset X\setminus U$, where $U=X\setminus\bigcup_{i=1}^s
\Supp(\Cal M/\fn_{i+r})$, $\J\not\subset\grad{0}$.
Since $\J\Cal L=0$ and $0$ is $G$-primary, $\Cal L=0$, as desired.
\qed

\begin{corollary}\label{G-primary-star.thm}
If $\Cal N$ is $G$-primary and
\[
\Cal N=\fn_1\cap\cdots\cap\fn_r
\]
is a minimal primary decomposition,
then $\fn_i^*=\Cal N$ for $i=1,\ldots,r$.
\end{corollary}

\proof Follows immediately from
Corollary~\ref{G-primary-no-emb.thm} and
Lemma~\ref{primary-star.thm}.
\qed

\begin{corollary}
If $\P$ is a $G$-prime $G$-ideal of $\O_X$,
then for any associated prime ideal $\fp$ of $\P$, $\fp^*=\P$.
\end{corollary}

\proof Follows immediately from
Corollary~\ref{G-primary-no-emb.thm} and
Corollary~\ref{G-prime-minimal}.
\qed

\paragraph\label{S_1-setting.par}
Assume that $X$ satisfies the $(S_1)$ condition (i.e., $\Cal O_X$ 
satisfies the $(S_1)$ condition).
Let
\begin{equation}\label{pd.eq}
0=\fq_1\cap\cdots\cap\fq_r
\end{equation}
be the minimal primary decomposition.
Set $X_i:=V(\fq_i)$, and $Y_i:=X_i\setminus\bigcup_{j\neq i}X_j$.
Define $G_{ij}=p_2^{-1}Y_i\cap a^{-1}Y_j$.
We say that $i\rightarrow j$ if $X_i^*\supset X_j$.

\begin{lemma}\label{G_ij.thm}
Let the notation be as above.
For $1\leq i,j\leq r$, the following are equivalent.
\begin{description}
\item[(i)] $i\rightarrow j$.
\item[(ii)] $G_{ij}\neq\emptyset$.
\end{description}
\end{lemma}

\proof {\bf (i)$\Rightarrow$(ii)} 
Since the closure of $G\times Y_i$ in $G\times X$ 
is $G\times X_i$, the scheme theoretic image of the action
$a|_{G\times Y_i}:G\times Y_i\rightarrow X$ is $X_i^*$.
Since $X_i^*\supset X_j$, $G\times Y_i$ intersects $a^{-1}(Y_j)$.
Namely, $G_{ij}\neq\emptyset$.

{\bf (ii)$\Rightarrow$(i)} 
Applying Corollary~\ref{S_1.thm}, 
the scheme theoretic image of $a|_{G_{ij}}:G_{ij}\rightarrow X_j$ is $X_j$.
This shows $X_i^*\supset X_j$.
\qed

\begin{corollary}
$\rightarrow$ is an equivalence relation on $\{1,\ldots,r\}$.
\end{corollary}

\proof Since $X_i^*\supset X_i$, $i\rightarrow i$.

Consider the isomorphism $\Theta:G\times X\rightarrow G\times X$ given
by $\Theta(g,x)=(g^{-1},gx)$.
Then $\Theta(G_{ij})=G_{ji}$.
Thus $i\rightarrow j$ if and only if $j\rightarrow i$.

Assume that $i\rightarrow j$ and $j\rightarrow k$.
Then $X_i^*=X_i^{**}\supset X_j^*\supset X_k$.
Hence $i\rightarrow k$.
\qed

\begin{lemma}\label{i->j.thm}
Assume that $\Cal M/\Cal N$ does not have an embedded prime ideal.
Let
\begin{equation}\label{n-primary.eq}
\Cal N=\fn_1\cap\cdots\cap \fn_r
\end{equation}
be a minimal primary decomposition.
Then
\[
\fn_i^*=\bigcap_{i\rightarrow j}\fn_j,
\]
where we say that $i\rightarrow j$ if $\uann(\Cal M/\fn_i)^*
\subset \uann(\Cal M/\fn_j)$.
\end{lemma}

\proof Replacing $\Cal M$ by $\Cal M/\Cal N$ and $\fn_j$ by $\fn_j/\Cal N$,
we may assume that $\Cal N=0$.
Replacing $X$ by $\Supp\Cal M:=V(\uann\Cal M)$, we may assume that
$\uann\Cal M=0$.
Set $\fq_j:=\uann(\Cal M/\fn_j)$ and $X_j:=V(\fq_j)$.
Note that (\ref{pd.eq}) is a minimal primary decomposition
by Lemma~\ref{M-ann.thm}, and $X$ satisfies
the $(S_1)$ condition.

Now we define $Y_j$ and $G_{ij}$ as in (\ref{S_1-setting.par}).
The definition of $i\rightarrow j$ is consistent with that in
(\ref{S_1-setting.par}).

Let $\rho_i:Y_i\rightarrow X$ be the inclusion, 
$\rho:\coprod_{i\rightarrow j}Y_j\rightarrow X$ be the inclusion,
$\varphi:\coprod_{i\rightarrow j} G_{ij}\rightarrow G\times X$ 
be the inclusion, 
$a_0:\coprod_{i\rightarrow j}G_{ij}\rightarrow \coprod_{i\rightarrow j}Y_j$ 
be the restriction of $a$, and
$p_0:\coprod_{i\rightarrow j}G_{ij}\rightarrow Y_i$ be the restriction
of $p_2$, so that the diagrams
\[
\xymatrix{
\coprod_{i\rightarrow j}G_{ij} \ar[r]^\varphi \ar[d]^{a_0} & 
G\times X \ar[d]^a \\
\coprod_{i\rightarrow j}Y_j \ar[r]^\rho &
X
}
\qquad
\xymatrix{
\coprod_{i\rightarrow j}G_{ij} \ar[r]^\varphi \ar[d]^{p_0} & 
G\times X \ar[d]^{p_2} \\
Y_i \ar[r]^{\rho_i} &
X
}
\]
are commutative.

Since $a(G\times Y_i)\subset \bigcup_{i\rightarrow j}X_j$ and 
$\coprod_{i\rightarrow j}Y_j$ is dense in $\bigcup_{i\rightarrow j}X_j$,
$\coprod_{i\rightarrow j}G_{ij}$ is dense in $G\times Y_i$ by 
(\ref{sim-flat.par}).
As $G\times Y_i$ is dense in $G\times X_i$, $\coprod_{i\rightarrow j}G_{ij}$
is also dense in $G\times X_i$.
Since $\Supp p_2^*(\Cal M/\fn_i)=G\times X_i$ and $p_2^*(\Cal M/\fn_i)$ 
satisfies the $(S_1)$ condition, 
$u:a_*p_2^*(\Cal M/\fn_i)\rightarrow
a_*\varphi_*\varphi^*p_2^*(\Cal M/\fn_i)$ is a monomorphism
by Lemma~\ref{S_1-flat.thm}.

Similarly, $u:\rho_*\rho^*\Cal M\rightarrow \rho_*(a_0)_*a_0^*\rho^*\Cal M$ 
is a monomorphism.

Since the diagram
\[
\xymatrix{
\rho_*\rho^*\Cal M \ar@{^{(}->}[r]^-u &
\rho_*(a_0)_*a_0^*\rho^*\Cal M \ar[r]^\cong &
a_*\varphi_*\varphi^*a^*\Cal M \ar[d]^\cong \\
\Cal M \ar[u]_u \ar[r]^\omega &
a_*p_2^*\Cal M \ar[r]^u \ar[d]^{\pi_{\fn_i}} &
a_*\varphi_*\varphi^*p_2^*\Cal M \ar[d]^\cong \\
 & a_*p_2^*(\Cal M/\fn_i) \ar[d]^u &
a_*\varphi_*p_0^*\rho_i^*\Cal M\ar[d]^\cong \\
 & a_*\varphi_*\varphi^*p_2^*(\Cal M/\fn_i)\ar[r]^\cong &
a_*\varphi_*p_0^*\rho_i^*(\Cal M/\fn_i)
}
\]
is commutative, 
\[
\fn_i^*=\Ker(\Cal M\xrightarrow\omega a_*p_2^*\Cal M
\xrightarrow{\pi_{\fn_i}} a_*p_2^*(\Cal M/\fn_i))
=\Ker(u:\Cal M\rightarrow \rho_*\rho^*\Cal M)=\bigcap_{i\rightarrow j}
\fn_j
\]
by Lemma~\ref{star-omega.thm}.
\qed

\begin{corollary}\label{all-i->j.thm}
Let $\Cal N$ be $G$-primary in {\rm Lemma~\ref{i->j.thm}}.
Then for $i,j\in\{1,\ldots,r\}$, $i\rightarrow j$.
\end{corollary}

\proof By Corollary~\ref{G-primary-star.thm} and Lemma~\ref{i->j.thm},
$\Cal N=\fn_i^*=\bigcap_{i\rightarrow j}\fn_j$.
Since the decomposition (\ref{n-primary.eq}) is irredundant, 
$i\rightarrow j$ holds for all $j\in\{1,\ldots,r\}$.
\qed

\begin{theorem}\label{main.thm}
Let 
\begin{equation}\label{N-M3.eq}
\Cal N=\Cal M_1\cap\cdots\cap \Cal M_s
\end{equation}
be a minimal $G$-primary decomposition, and 
let
\begin{equation}\label{M-m.eq}
\Cal M_l=\fm_{l,1}\cap\cdots\cap\fm_{l,r_l}
\end{equation}
be a minimal primary decomposition.
Then
\[
\Cal N=\bigcap_{l=1}^s(\fm_{l,1}\cap\cdots\cap\fm_{l,r_l})
\]
is a minimal primary decomposition.
\end{theorem}

\proof Note that $\sqrt{\fm_{l,j}:\Cal M}$ are distinct.
Indeed, we have
\[
\sqrt{\fm_{l,j}:\Cal M}^{\,*}
=
\sqrt[G]{(\fm_{l,j}:\Cal M)^*}
=\sqrt[G]{\fm_{l,j}^*:\Cal M}
=\sqrt[G]{\Cal M_l:\Cal M}.
\]
Since (\ref{N-M3.eq}) is minimal, $\sqrt{\fm_{l,j}:\Cal M}$ is
different if $l$ is different.
On the other hand, if $l$ is the same and $j$ is different, then
$\sqrt{\fm_{l,j}:\Cal M}$ is different, since (\ref{M-m.eq}) is minimal.

So it suffices to prove that each $\sqrt{\fm_{l,j}:\Cal M}$ is an
associated prime ideal of $\Cal N$.
By Lemma~\ref{module-ideal-G-primary.thm}, $\Cal M_l:\Cal M$ is
$G$-primary, and $\sqrt[G]{\Cal M_l:\Cal M}$ is $G$-prime.
So neither $\Cal M_l:\Cal M$ nor $\sqrt[G]{\Cal M_l:\Cal M}$ has an
embedded prime ideal by Corollary~\ref{G-primary-no-emb.thm}.
So
\[
\Ass(\Cal O_X/(\Cal M_l:\Cal M))=\Ass(\Cal O_X/(\sqrt[G]{\Cal M_l:\Cal M}))=
\Ass(\Cal O_X/\sqrt{\Cal M_l:\Cal M}).
\]
By Lemma~\ref{M-ann.thm}, 
\[
\Cal M_l:\Cal M=\bigcap_j\fm_{l,j}:\Cal M
\]
is a minimal primary decomposition.
So $\sqrt[G]{\Cal M_l:\Cal M}\in\Ass_G\Cal M/\Cal N$ and
$\sqrt{\fm_{l,j}:\Cal M}\in\Ass(\Cal O_X/\sqrt[G]{\Cal M_l:\Cal M})$.
By Lemma~\ref{Ass-Ass.thm}, 
$\sqrt{\fm_{l,j}:\Cal M}\in\Ass(\Cal M/\Cal N)$,
as desired.
\qed

\begin{corollary}\label{ass-star.thm}
A prime ideal $\fp$ of $\Cal O_X$ is an associated prime ideal of some 
coherent $(G,\Cal O_X)$-module $\Cal M$ if and only if $\fp$ is a
minimal prime ideal of $\fp^*$.
\end{corollary}

\proof The \lq if' part is trivial.
We prove the converse.
Take a minimal $G$-primary decomposition (\ref{N-M3.eq}) and
minimal primary decompositions (\ref{M-m.eq}).
We may assume that $\fp=\sqrt{\fm_{1,1}:\M}$.
Then $\fp^*=\grad{\M_1:\M}$.
As $\M/\M_1$ does not have an embedded prime ideal, 
$\fm_{1,1}:\M$ is a primary component of $\M_1:\M$ corresponding to
a minimal prime ideal, and hence
$\fp$ is a minimal prime ideal of $\M_1:\M$.
Since $\M_1:\M$ and $\fp^*=\grad{\M_1:\M}$ have the same radical, 
we are done.
\qed

\begin{corollary}\label{Ass-Ass2.thm}
Let {\rm(\ref{N-M3.eq})} be a minimal $G$-primary decomposition.
Then
\[
\Ass(\Cal M/\Cal N)
=
\coprod_{l=1}^s\Ass(\Cal M/\Cal M_l)
=
\coprod_{\Cal P\in\Ass_G(\Cal M/\Cal N)}\Ass(\Cal O_X/\Cal P)
\]
and
\[
\Ass_G(\Cal M/\Cal N)=\{\fp^*\mid \fp\in\Ass(\Cal M/\Cal N)\}.
\]
\end{corollary}

\proof Follows immediately by Theorem~\ref{main.thm}.
\qed

\begin{corollary}\label{Min-Min.thm}
We have
\[
\Min(\Cal M/\Cal N)=\coprod_{\Cal P\in\Min_G(\Cal M/\Cal N)}
\Ass(\Cal O_X/\Cal P)
\]
and
\[
\Min_G(\Cal M/\Cal N)
=
\{\fp^*\mid \fp\in\Min(\Cal M/\Cal N)\}
\]
\end{corollary}

\proof
Assume that $\fp\in\Min(\Cal M/\Cal N)$ and $\fp^*\notin
\Min_G(\Cal M/\Cal N)$.
Then there exists some $\Cal P\in\Min_G(\Cal M/\Cal N)$ such that
$\Cal P\subsetneq \fp^*$.
If
\[
\sqrt{\Cal P}=\fp_1\cap\cdots\cap\fp_s
\]
is a minimal prime decomposition, then each $\fp_i$ is an element of 
$\Ass(\Cal M/\Cal N)$ by Corollary~\ref{Ass-Ass2.thm}.
Since
\[
\fp_1\cap\cdots\cap\fp_s=\sqrt{\Cal P}\subset\sqrt{\fp^*}\subset
\sqrt{\fp}=\fp,
\]
there exists some $i$ such that $\fp_i\subset\fp$.
Since $\fp_i^*=\Cal P\neq \fp^*$, we have $\fp_i\subsetneq \fp$.
This contradicts the minimality of $\fp$.
So $\fp\in\Min(\Cal M/\Cal N)$ implies $\fp^*\in\Min_G(\Cal M/\Cal N)$.
By Corollary~\ref{ass-star.thm}, $\fp\in\Min(\Cal M/\Cal N)$ implies
that $\fp$ is an associated prime ideal of $\fp^*$.
So the $\subset$ direction of the first equation 
and $\supset$ direction of the second equation
have been proved.

Conversely, assume that $\Cal P\in\Min_G(\Cal M/\Cal N)$,
and $\fp\in\Ass(\Cal O_X/\Cal P)$.
By Lemma~\ref{Ass-Ass.thm}, $\fp\in\Ass(\Cal M/\Cal N)$.
Assume that $\fp$ is not minimal.
Then there exists some $\fp'\in\Min(\Cal M/\Cal N)$ such that
$\fp'\subsetneq \fp$.
Then $(\fp')^*\subset \fp^*=\Cal P$.
By the minimality of $\Cal P$, 
$(\fp')^*=\fp^*=\Cal P$.
So $\fp'\in\Ass(\Cal O_X/\Cal P)$ by Corollary~\ref{ass-star.thm}.
As $\Cal O_X/\Cal P$ does not have an embedded prime ideal
by Corollary~\ref{G-primary-no-emb.thm}, this contradicts $\fp'\subsetneq
\fp$.
So $\fp$ must be minimal.
This proves the $\supset$ direction of the first equation.
As $\Cal P=\fp^*$ with $\fp\in\Min(\Cal M/\Cal N)$, 
the $\subset$ direction of the second equation has also been proved.
\qed

\begin{corollary}\label{no-emb-no-emb.thm}
We have $\Ass(\M/\N)=\Min(\M/\N)$ if and only if $\Ass_G(\M/\N)
=\Min_G(\M/\N)$.
\end{corollary}

\proof Obvious by Corollary~\ref{Ass-Ass2.thm} and
Corollary~\ref{Min-Min.thm}.
\qed

\begin{corollary}\label{G-radical-no-emb.thm}
A $G$-radical $G$-ideal does not have an embedded prime ideal.
\end{corollary}

\proof Obvious by Corollary~\ref{no-emb-no-emb.thm}.
\qed

\begin{example}
Let $k=\Bbb Z/2\Bbb Z$ and $S=\Spec k$.
Let $G=\mu_6=\Spec H$,
the group scheme of the 6th roots of unity,
where $H=k[t]/(t^6-1)$.
Let $A=k[x,y,z]$ and $I=((x-y)^3,(x-y)(x-z))$.
$G$ acts on $A$ via the coaction $\omega:A\rightarrow A\otimes H$
given by $\omega(x)=x\otimes 1$, $\omega(y)=y\otimes t$, and $\omega(z)=
z\otimes t$.
Letting $J=(x-y)$ and $K=((x-y)^3,x-z)$, $I=J\cap K$ is a primary 
decomposition.
By Lemma~\ref{star-omega.thm}, it is easy to see that
$L^*$ is the kernel of the map
\[
A \xrightarrow{\omega}A\otimes H \xrightarrow{p\otimes 1} A/L\otimes H,
\]
where $p:A\rightarrow A/L$ is the projection.

Using Macaulay2 \cite{M2}, it is easy to compute
\[
J^*=(x^6-y^6)=(x^2-y^2)\cap(x^4+x^2y^2+y^4)
\]
and 
\begin{multline*}
K^*=(y^3+y^2z+yz^2+z^3,x^6-z^6)
=\\
(x^2-z^2,y^3+y^2z+yz^2+z^3)\cap
(y^3+y^2z+yz^2+z^3,x^3+x^2z^2+z^4)
\end{multline*}
are minimal primary decompositions.

As $J$ and $K$ are primary, $J^*$ and $K^*$ are $G$-primary.
It is easy to check that the $G$-primary decomposition
$I^*=J^*\cap K^*$ is irredundant, and minimal in fact.
So by Theorem~\ref{main.thm},
\begin{multline*}
I^*=
(x^2-y^2)\cap(x^4+x^2y^2+y^4)
\cap\\
(x^2-z^2,y^3+y^2z+yz^2+z^3)\cap
(y^3+y^2z+yz^2+z^3,x^3+x^2z^2+z^4)
\end{multline*}
is a minimal primary decomposition.
\end{example}

\begin{lemma}\label{nowhere.thm}
If $X$ is $G$-integral and $\I$ a nonzero $G$-ideal, then $V(\I)$ is
nowhere dense in $X$.
\end{lemma}

\proof Assume the contrary.
Then there exists some minimal prime ideal
$\fp$ of $0$ such that $\I\subset\fp$.
Then $\I=\I^*\subset \fp^*=0$, since $0$ is $G$-prime.
This is a contradiction.
\qed

\begin{lemma}\label{G-integral-projective-dense.thm}
Let $X$ be $G$-integral, and $\M$ a coherent
$(G,\O_X)$-module.
Then there exists some $r$ and some dense $G$-stable open subset $U$ of $X$
such that for $x\in X$, it holds $x\in U$ if and only if $\M_x\cong \O_{X,x}^r
$.
In this case, $\M|_U$ is locally free of rank $r$.
\end{lemma}

\proof Let $r$ be the smallest integer such that $\uFitt_r\M\neq0$.
Then $r\geq 0$, and letting $U:=X\setminus V(\uFitt_r\M)$, $U$ is 
$G$-stable open, and for $x\in X$, it holds $x\in U$ if and only if
$\M_x\cong \O_{X,x}^r$ by
\cite[Proposition~20.8]{Eisenbud}.

$U$ is dense, since $V(\uFitt_r\M)$ is nowhere dense by 
Lemma~\ref{nowhere.thm} and closed.
\qed

\begin{lemma}\label{G-reduced-free.thm}
Let $X$ be $G$-reduced, and $\M$ a coherent $(G,\O_X)$-module.
Then $U:=\{x\in X\mid \M_x \text{ is projective}\}$ is a dense $G$-stable
open subset of $X$.
\end{lemma}

\proof It is easy to see 
that $U=\bigcup_{r\geq 0}(X\setminus (V(\uFitt_r\M)\cup\Supp
\uFitt_{r-1}\M))$ is $G$-stable open, and if $x\in U$, then $\M_x$ is 
projective.
Conversely, if $\M_x$ is projective of rank $r$, then $\Fitt_r(\M_x)=\O_{X,x}
$ and $\Fitt_{r-1}(\M_x)=0$ by \cite[Proposition~20.8]{Eisenbud}, 
and hence $x\in U$.

It remains to show that $U$ is dense.
Let $\fp$ be any associated (or equivalently, minimal, by
Corollary~\ref{G-radical-no-emb.thm}) prime ideal of $0$.
We need to show that the generic point $\xi$ of $V(\fp)$ is in $U$.

Let $0=\P_1\cap \cdots\cap \P_r$ be a minimal $G$-prime decomposition.
Then we may assume that $\fp$ is a minimal prime ideal of $\P_1$.
As $X$ is $G$-reduced, $\fp\not\supset\P_i$ for $i\geq 2$.
Hence $\xi\in Y:=X\setminus(
\bigcup_{i\geq 2}V(\P_i))$.
As $Y$ is a 
non-empty
$G$-stable open subscheme of $V(\P_1)$,
it is $G$-integral by Lemma~\ref{G-integral-open.thm},
and $\xi$ is a generic point
of an irreducible component of $Y$.
Since $Y\cap U$ is dense in $Y$ by
Lemma~\ref{G-integral-projective-dense.thm} and its proof, 
$\xi\in U$, as desired.
\qed

\begin{corollary}
Let $X$ be $G$-reduced, and $\L$ a quasi-coherent 
$(G,\O_X)$-module.
Then for the generic point $\xi$ of an irreducible component of $X$, 
$\L_\xi$ is $\O_{X,\xi}$-flat.
\end{corollary}

\proof Since $\L$ is a filtered inductive limit $\indlim \M_\lambda$
of its coherent $(G,\O_X)$-submodules $\M_\lambda$ by 
(\ref{locally-noetherian.par}) 
and
$(\M_\lambda)_{\xi}$ is a free module by Lemma~\ref{G-reduced-free.thm},
$\L_\xi$ is $\O_{X,\xi}$-flat.
\qed

\begin{corollary}
Let $X$ be $G$-reduced, and $f:V\rightarrow X$ an
affine $G$-morphism.
Let $v\in V$, and assume that $f(v)$ is a generic point of $X$, then 
$f$ is flat at $v$.
\end{corollary}

\proof This is because $(f_*\O_V)_{f(v)}$ is $\O_{X,f(v)}$-flat, and
$\O_{V,v}$ is a localization of $(f_*\O_V)_{f(v)}$.
\qed

\begin{lemma}\label{smooth.thm}
Let $p_2:G\times X\rightarrow X$ have regular fibers.
Then 
\begin{description}
\item[(i)] If $\fa$ is a radical quasi-coherent ideal of $\Cal O_X$,
then $\fa^*$ is also radical.
\item[(ii)] Any $G$-radical $G$-ideal of $\Cal O_X$ is radical.
\end{description}
\end{lemma}

\proof Clearly, {\bf (ii)} follows from {\bf (i)}.
We prove {\bf (i)}.

Set $Y:=V(\fa)$.
Then $Y$ is reduced.
By assumption, $G\times Y$ is reduced.
So the scheme theoretic image $Y^*$ of the action $a:G\times Y\rightarrow X$
is also reduced.
Since $Y^*=V(\fa^*)$, we have that $\fa^*$ is radical.
\qed

\begin{corollary}
Let $p_2:G\times X\rightarrow X$ have regular fibers.
If $\I$ is a $G$-ideal of $\O_X$, then $\sqrt\I=\grad\I$ 
is a $G$-radical $G$-ideal.
\end{corollary}

\proof Note that $\grad\I$ is $G$-radical, and hence is radical by the lemma.
Hence
\[
\sqrt\I=\sqrt{\grad\I}=\grad\I
\]
is a $G$-radical $G$-ideal.
\qed

\begin{lemma}\label{primary2.thm}
Let $p_2:G\times X\rightarrow X$ have connected fibers.
Then
\begin{description}
\item[(i)] If $\fn$ is a primary submodule of $\M$,
then $\fn^*$ is also primary.
\item[(ii)] Any $G$-primary $G$-submodule of $\M$ is a primary submodule.
\end{description}
\end{lemma}

\proof Since {\bf (ii)} follows from {\bf (i)} 
and Corollary~\ref{G-primary-star.thm},
we only need to prove {\bf (i)}.

By assumption, $\M/\fn$ is primary, and hence it satisfies $(S_1)$.
As $p_2:G\times X\rightarrow X$ is flat with Cohen--Macaulay fibers
(see e.g., \cite[(31.14)]{ETI}), 
$p_2^*(\M/\fn)$ also satisfies $(S_1)$
by \cite[(6.4.1)]{EGA-IV}.
As $p_2$ is faithfully flat, $\Supp p_2^*(\M/\fn)=p_2^{-1}(\Supp(\M/\fn))$,
which is irreducible by the irreducibility of $\Supp(\M/\fn)$, and
the assumption that $p_2$ has connected (or equivalently, geometrically 
irreducible) fibers.
Thus $p_2^*(\M/\fn)$ is a primary module.
By Lemma~\ref{star-omega.thm}, there is a monomorphism
$\M/\fn^*\hookrightarrow a_*p_2^*(\M/\fn)$.
By Lemma~\ref{primary-sim.thm}, $\fn^*$ is a primary submodule, as desired.
\qed

\begin{corollary}\label{primary3.thm}
Let $p_2:G\times X\rightarrow X$ have connected fibers.
If $\M$ is a coherent $\O_X$-module and $\N$ is a
coherent $(G,\O_X)$-submodule, then
a minimal $G$-primary decomposition of $\N$ is also a minimal primary
decomposition of $\N$.
\end{corollary}

\proof Follows from Lemma~\ref{primary2.thm} and
Theorem~\ref{main.thm}.
\qed

\begin{corollary}\label{regular-and-connected.thm}
Let $p_2:G\times X\rightarrow X$ have regular and 
connected fibers.
Then
\begin{description}
\item[(i)] If $\fp$ is a prime ideal of $\cO_X$, then $\fp^*$ is also a
prime ideal.
\item[(ii)] Any $G$-prime $G$-ideal of $\cO_X$ is prime.
\item[(iii)] For a coherent $(G,\O_X)$-module $\M$ of $\O_X$, 
$\Ass_G(\M)=\Ass(\M)$ and
$\Min_G(\M)=\Min(\M)$.
In particular, any associated prime ideal of a coherent 
$(G,\O_X)$-module is a $G$-prime $G$-ideal.
\end{description}
\end{corollary}

\proof {\bf (i)} and {\bf (ii)} 
follow immediately from Lemma~\ref{smooth.thm} and
Lemma~\ref{primary2.thm}.
{\bf (iii)} follows from {\bf (ii)}, Corollary~\ref{Ass-Ass2.thm}, and
Corollary~\ref{Min-Min.thm}.
\qed

\begin{example}
Let $S=\Spec \Bbb Z$, $G=\Bbb G_m^n$, and $X=\Spec A$.
Then $A$ is a noetherian $\Bbb Z^n$-graded ring.
For a prime ideal $P$ of $A$, $P^*$, the ideal generated by all the
homogeneous elements of $P$ is a prime ideal by
Corollary~\ref{regular-and-connected.thm}.
An associated ideal of a homogeneous ideal is again homogeneous.
These facts are well-known, and checked directly.
\end{example}

\begin{lemma}\label{locally-constant.thm}
Let $S$ be locally noetherian and $G$ be flat of finite type over $S$.
Let $\pi$ be the structure map $G\rightarrow S$.
Then $h(s):=\dim(\pi^{-1}(s))$ is a locally constant function on $S$.
\end{lemma}

\proof We may assume that $S=\Spec A$ is affine and reduced.
Then it suffices to show that $h$ is constant on each irreducible component.
So we may further assume that $S$ is irreducible.
Let $\sigma$ be the generic point of $S$.
The fiber $\pi^{-1}(\sigma)$ is a finite-type
group scheme over the field $\kappa(\sigma)$, and is equidimensional
of dimension $h(\sigma)$.

Let $s\in S$.
Note that
$\pi$ is an open map (follows easily from \cite[(1.10.4)]{EGA-IV}).
By \cite[(14.2.4)]{EGA-IV}, each irreducible component of $\pi^{-1}(s)$ is
$h(\sigma)$-dimensional.
Hence, $h(s)=h(\sigma)$, and this value is independent of $s$.
\qed

\begin{proposition}\label{constant.thm}
Assume that $\Min_G(0)$ is a singleton, where $0$ is the zero ideal of $\O_X$.
Then the dimension of the fiber of $p_2:G\times X\rightarrow X$ is constant.
In particular, if $X$ is $G$-primary, then
the dimension of the fiber of $p_2:G\times X\rightarrow X$ is constant.
\end{proposition}

\proof 
Set $h(x)=\dim p_2^{-1}(x)$.
We want to prove that $h$ is constant.
By Lemma~\ref{locally-constant.thm}, $h$ is locally constant.
Let $0=\Q_1\cap\Q_2\cap\cdots\cap\Q_r$ be a minimal $G$-primary decomposition
where $\grad{\Q_1}$ is a minimal $G$-prime $G$-ideal.
Then $Y=V(\Q_2\cap\cdots\cap\Q_r)$ is nowhere dense in $X$.
By the local constantness, replacing $X$ by $X\setminus Y$, we may assume
that $X$ is $G$-primary.
Let $0=\fq_1\cap \cdots \cap
\fq_s$ be the minimal primary decomposition of $0$ 
in $\O_X$.
Set $X_i:=V(\fq_i)$, and let $\xi_i$ be the generic point of $X_i$.
Let $\pi:X\rightarrow S$ be the structure map.
Let $Y_i:=X_i\setminus\bigcup_{j\neq i}X_j$.

It suffices to show that for $1\leq i,j\leq s$, $h(\xi_i)=h(\xi_j)$.
Note that $G_{ij}=a^{-1}Y_j\cap p_2^{-1}Y_i$ is non-empty
by Corollary~\ref{all-i->j.thm} and 
Lemma~\ref{G_ij.thm}.
Let $\gamma$ be the generic point of an irreducible component of $G_{ij}$.
By flatness, $a(\gamma)=\xi_j$ and $p_{2}(\gamma)=\xi_i$.
Hence $\sigma:=\pi(\xi_j)=\pi a(\gamma)=\pi p_2(\gamma)=\pi(\xi_i)$.
Since $\{\xi_j\}\rightarrow \{\sigma\}$ associates with a field extension, 
and is faithfully flat quasi-compact, and $G\times\{\xi_j\}$ is of finite type
over $\{\xi_j\}$, we have that $G\times\{\sigma\}$ is 
of finite type
over $\{\sigma\}$
by \cite[(2.7.1)]{EGA-IV}.
It is easy to see that $G\times\{\sigma\}$ and $G\times\{\xi_j\}$ have
the same dimension, $h(\xi_j)$.
By the same reason, 
$G\times\{\sigma\}$ is $h(\xi_i)$-dimensional.
Hence $h(\xi_i)=h(\xi_j)$.
\qed

\paragraph\label{Ratliff.par}
A commutative ring $A$ is said to be {\em Hilbert} if 
any prime ideal $P$ of $A$ equals the intersection $\bigcap_{\fm}
\fm$, where the intersection is taken over the maximal ideals containing $P$.
We say that $A$ satisfies the {\em first chain condition} (FCC for short)
if each maximal chain of prime ideals in $A$ has the length equals to $\dim A$.
We say that $A$ is {\em Ratliff} if $A$ is noetherian, universally catenary, 
Hilbert, and for
any minimal prime ideal $P$ of $A$, $A/P$ satisfies the FCC.

\begin{lemma}
A noetherian ring $A$ is Ratliff if and only if $A\red$ is.
\end{lemma}

\proof By assumption, both $A$ and $A\red$ are noetherian.
Almost by definition, $A$ is Hilbert if and only if $A\red$ is.
A finite-type algebra over $A\red$ is of finite type over $A$.
So plainly, if $A$ is universally catenary, then $A\red$ is universally
catenary.
Let $B$ be a finite type algebra over $A$.
Then $B$ is catenary if and only if $B\red$ is, and $B\red$ is a
finite-type algebra over $A\red$.
So if $A\red$ is universally catenary, then so is $A$.
Note that $\Min(A\red)=\{PA\red\mid P\in\Min(A)\}$, and
$A/P\cong A\red/PA\red$.
So $A/P$ satisfies the FCC for any $P\in\Min(A)$ if and only if
$A\red/Q$ satisfies the FCC for any $Q\in \Min(A\red)$.
\qed

\paragraph
An artinian ring is Ratliff.
A one-dimensional noetherian domain with infinitely many prime ideals is
Ratliff.
For example, $\Bbb Z$ is Ratliff.

\begin{lemma}\label{FCC.thm}
Let $A$ be a catenary noetherian ring such that for each
minimal prime ideal $P$ of $A$, $A/P$ satisfies the FCC.
Then for any prime ideal $Q$ of $A$, $A/Q$ satisfies the FCC.
\end{lemma}

\proof Easy.
\qed

\begin{corollary} A homomorphic image of a Ratliff ring is Ratliff.
\end{corollary}

\proof Follows from Lemma~\ref{FCC.thm}.
\qed

\begin{lemma}\label{Ratliff-of-finite-type.thm}
If $A$ is Ratliff and $A\rightarrow B$ is of finite type, then $B$ is 
Ratliff.
\end{lemma}

\proof Note that $B$ is noetherian by Hilbert's basis theorem.
$B$ is universally catenary, since a finite-type algebra $C$ over $B$
is also of finite type over $A$.
It is well-known that a finite-type algebra over a Hilbert ring is
Hilbert \cite[Chapter~6, Theorem~1]{Northcott}.
It remains to show that for any minimal prime 
ideal $P$, $\bar B:=B/P$ satisfies the 
FCC.
Since we know that $\bar B$ is catenary, 
it suffices to show that $\dim \bar B_{\fm}$ is independent of the choice 
of  a maximal ideal $\fm$ of $\bar B$.
Let $\fp:=P\cap A$, and set $\bar A:=A/\fp$.
Note that $\bar A$ is Ratliff.
Note also that $\fn:=\fm\cap \bar A$ is a maximal ideal of $\bar A$ 
\cite[Chapter~6, Theorem~2]{Northcott},
and $\dim \bar A_{\fn}=\dim \bar A$ is independent of $\fn$.
Clearly, $\kappa(\fm)$ is an algebraic extension of $\kappa(\fn)$.
By the dimension formula \cite[(5.6.1)]{EGA-IV},
$\dim \bar B_{\fm}=\dim \bar A+\tdeg_{\bar A}\bar B$ is independent of 
$\fm$.
\qed

\paragraph
We say that a scheme $Y$ is Ratliff if $Y$ has a finite affine open covering
$(U_i)$ such that each $U_i$ is isomorphic to the prime spectrum of a 
Ratliff ring.
$Y$ is Ratliff if and only if $Y$ is quasi-compact, and each point
of $Y$ has an affine Ratliff open neighborhood.
A Ratliff scheme $Y$ has a finite dimension.

If $Y$ is Ratliff and $Z\rightarrow Y$ is of finite type, then $Z$ is
Ratliff.

A noetherian scheme $Y$ is said to be {\em equidimensional}, if 
the irreducible components of $Y$ have the same dimension.
If $Y$ is Ratliff and equidimensional, then $\dim Y=\dim\O_{Y,y}$ for
any closed point $y$ of $Y$.
Let $Z$ be a closed subset of a Ratliff scheme $Y$.
Let $W$ be the set of closed points of $Z$.
Then the closure of $W$ is $Z$.
It follows that any 
non-empty open subscheme $U$ of a Ratliff scheme $Y$ 
contains a closed point of $Y$.
In particular, if $Y$ is Ratliff equidimensional 
and $U$ is a non-empty open subset, then $\dim U=\dim Y$.
In particular, for any point $y$ of $Y$, $\dim_y Y=\dim Y$,
where $\dim_y Y=\inf\{\dim U\mid U \text{ is an open neighborhood of $y$ 
in $Y$}\}$, see \cite[(0.14.1.2)]{EGA-IV}.

Let $f:Z\rightarrow Y$ be of finite type and dominating, $Y$ be Ratliff, and
$Z$ and $Y$ be irreducible.
Then $\dim Z=\dim Y+\tdeg_{\kappa(\eta)}\kappa(\zeta)$, 
where $\eta$ and $\zeta$ are respectively the generic points of $Y$ and $Z$.

\begin{proposition}\label{equidimensional.thm}
Assume that $X$ is Ratliff, and assume that $\Min_G(0)$ is a singleton, 
where $0$ is the zero ideal of $\O_X$.
Then $X$ is equidimensional.
\end{proposition} 

\proof Discarding a nowhere dense $G$-stable closed subscheme form $X$,
we may assume that $X$ is $G$-primary.

Let $X_i$, $Y_i$, and $G_{ij}$ be as in the proof of 
Proposition~\ref{constant.thm}.
Let $g$ be a closed point of $G_{ij}$.
Then $p_2(g)$ is a closed point of $Y_i$, and we have
$\dim \O_{G,g}=\dim Y_i+h$, where $h$ is the dimension of the fibers of
$p_2:G\times X\rightarrow X$, see Proposition~\ref{constant.thm}.
Since $a(g)$ is a closed point of $Y_j$, 
$\dim \O_{G,g}=\dim Y_j+h$ (note that the dimension of the fibers of
$a$ also have the constant value $h$).
Hence, $\dim X_j=\dim Y_j=\dim Y_i=\dim X_i$, as desired.
\qed

\begin{example}
Let $R$ be a DVR, and $t$ a prime element of $R$, and $K=R[t^{-1}]$ 
the field of fractions of $R$.
Set $S=\Spec R$.
Let $\tilde G= S_0 \coprod S_1 $ 
be the constant group $\Bbb Z/2\Bbb Z$ over $S$,
where $S=S_0=S_1$, and $S_0$ corresponds to the unit element, and 
$S_1$ corresponds to the other element of $\Bbb Z/2\Bbb Z$.
Then $G=S_0\coprod \Spec K$ is a flat of finite type group scheme over $S$.
Letting $X=G$, the left regular action, $X$ is $G$-primary.
But $X$ is not equidimensional.

This example also shows that the following statement is {\em false}.
Let $\P$ and $\Q$ be $G$-prime $G$-ideals of $\O_X$ such that $\Q\supset \P$.
Then for any minimal prime ideal $\fp$ of $\P$, there exists some minimal prime
ideal
$\fq$ of $\Q$ such that $\fq\supset \fp$.
Indeed, let $\P=0$, 
and $\Q$ be the defining ideal of the closed point of $S_0$.
Then the component $\Spec K\subset S_1$ of $V(\P)$ does not 
contain a component of
$V(\Q)$.
\end{example}

\begin{example}
Let $(R,t)$ be a DVR of mixed characteristic $p$.
Then $\mu_{p,R}=\Spec R[x]/(x^p-1)$ is a local scheme 
(i.e., the prime spectrum of a local ring, see \cite[(2.4.1)]{EGA-I}) with
two or more irreducible components.
\end{example}

\begin{example}
A primary ideal $\fq$ of $\O_X$ which is a primary component of some 
$G$-ideal  of $\O_X$, but not a primary component of $\fq^*$.
\end{example}

\proof[Construction.] Let $S=\Spec k$ with $k$ a field.
Let $G=\Bbb G_m^2$, and $X=\Spec A$, $A=k[x,y]$.
$G$ acts on $A$ with $\deg x=(1,0)$ and $\deg y=(0,1)$.
Let $\fq=(x^4,yx^3,y^2x^2+y^3x,y^4)$.
It is a primary component of the homogeneous ideal
$I=(x^4,yx^3)=(x^3)\cap\fq$.
But $\fq$ is not a primary component of $\fq^*=(x^4,yx^3,y^3x^2,y^4)$.
\qed

\begin{lemma}
Let $\N$ be a $(G,\O_X)$-submodule of $\M$.
If $\fm$ is the primary component of $\N$ corresponding to a 
minimal prime ideal
of $\N$, then $\fm$ is the primary component of $\fm^*$ corresponding
to a minimal prime ideal of $\fm^*$.
\end{lemma}

\proof
Let (\ref{N-M3.eq}) be a minimal $G$-primary decomposition of $\N$, and
(\ref{M-m.eq}) be a minimal primary decomposition of $\M_l$.
By the uniqueness of the primary component for minimal prime ideals,
we may assume that $\fm_{1,1}=\fm$.
Then $\fm$ is a primary component of $\M_1=\fm^*$ corresponding to a minimal
prime ideal by Corollary~\ref{G-primary-star.thm} and
Corollary~\ref{G-primary-no-emb.thm}.
\qed

\begin{example}
Even if $\fm$ is a maximal ideal of $\O_X$, $\fm^*$ may not be a 
$G$-maximal $G$-ideal.
Let $S=\Spec k$ with $k$ a field, $G=\Bbb G_m$, and $X=\Bbb A^1$ on
which $G$ acts by multiplication.
For any maximal ideal $\fm$ of $k[X]=k[t]$ not corresponding to the origin, 
$\fm^*=0$ is not $G$-maximal.
\end{example}

\paragraph\label{PT.par}
Let us consider the case that $S=\Spec k$, where $k$ is an
algebraically closed field.
Let $G$ be a linear (smooth) algebraic group over $k$.
Let $G^\circ$ denote the identity component of $G$.
Then there exists some $h_0,h_1,\ldots,h_m\in G(k)$ such that $h_0=e$ is the
unit element of $G$, and $G=h_0G^\circ\coprod h_1G^\circ\coprod
\cdots\coprod h_mG^\circ$.
Let $\N$ be a coherent $(G,\O_X)$-submodule of $\M$.
Let $\fp\in\Ass(\M/\N)$.
Then $\fp$ is $G^\circ$-stable by Corollary~\ref{regular-and-connected.thm}.
Let $Y_{\fp,0}=V(\fp)$, and $Y_{\fp,i}=h_iY_{\fp,0}$ for $i=1,\ldots ,m$.
Then 
\begin{equation}\label{pt.eq}
Y_\fp=\bigcup_{i=0}^m Y_{\fp,i}
\end{equation}
with the reduced structure is $G$-stable.
So $Y_{\fp,i}^*=Y_\fp$ for each $i$, as can be seen easily.
So $\Ass_G(\M/\N)=\{\I(Y_\fp)\mid \fp\in\Ass(\M/\N)\}$
by Corollary~\ref{Ass-Ass2.thm}, where 
$\I(Y_\fp)$ is the defining ideal of $Y_\fp$.

We can show that for a coherent $(G,\O_X)$-module $\L$, 
$\L$ is $\I(Y_\fp)$-$G$-primary if and only if $\Ass \L=\{\I(Y_{\fp,i})\mid
i=0,\ldots,m\}$.
To verify this, note that $Y_{\fp,i}$ in (\ref{pt.eq}) are isomorphic 
one another, and (\ref{pt.eq}) is an irreducible decomposition without
embedded component (but there may be some redundancy).
The \lq only if' part follows from the fact that 
$\L$ does not have an embedded prime ideal 
by Corollary~\ref{G-primary-no-emb.thm},
and $V(\sqrt{0:\L})=Y_\fp$.
The if part follows from the fact that $\I(Y_{\fp,i})^*=\I(Y_\fp)$ for
$i=0,\ldots,m$.

In particular, the existence of $G$-primary decomposition of $\N$
gives another proof of 
\cite[Thoerem~4.18]{PT} for the case that $X$ is quasi-compact.

\section{Matijevic--Roberts type theorem}

In this section, $S$, $G$, $X$, and $\M$ are as in the last section.
As in the last
section, we assume that $p_2:G\times X\rightarrow X$ is
of finite type.

\paragraph A {\rm flat}
homomorphism of noetherian rings $\varphi:A\rightarrow B$ is said to be
l.c.i.\ (local complete intersection) 
(resp.\ regular) if for any prime ideal $P$ of $A$, the fiber ring
$B_P/PB_P$ is l.c.i.\ (resp.\ geometrically regular over the field
$A_P/PA_P$).
By the openness of l.c.i.\ locus \cite[(I.2.12.4)]{Hashimoto} 
(resp.\ smooth locus \cite[(6.8.7)]{EGA-IV}) for a finite-type
morphism, a flat local homomorphism
essentially of finite type $(A,\fm)\rightarrow (B,\fn)$ is l.c.i.\ 
(resp.\ regular) if and only if the closed fiber $B/\fm B$ is a 
complete intersection (resp.\ geometrically regular).

\begin{theorem}\label{M-R-orig.thm}
Let $y$ be a point of $X$, and $Y$ the integral closed subscheme of $X$ whose
generic point is $y$.
Let $\eta$ be the generic point of an irreducible component of $Y^*$.
Then there are a noetherian local ring $A$ and flat l.c.i.\ 
local homomorphisms
essentially of finite type
$\varphi: \Cal O_{X,y}\rightarrow A$ and
$\psi: \Cal O_{X,\eta}\rightarrow A$,
such that $\dim A=\dim\Cal O_{X,y}$ and that $A\otimes_{\O_{X,\eta}}\L_\eta
\cong A\otimes_{\O_{X,y}}\L_y$ for any quasi-coherent $(G,\O_X)$-module
$\L$ of $X$.
If, moreover, $p_2:G\times X\rightarrow X$ is smooth, then $\varphi$ and
$\psi$ can be taken to be regular.
\end{theorem}

\proof The action $a:G\times Y\rightarrow Y^*$ is dominating, so 
there is a point $z\in G\times Y$ such that $a(z)=\eta$.
Let $\zeta$ be the generic point of an irreducible component of $G\times Y$
containing $z$.
Then $a(\zeta)$ is a generalization of $\eta$, and hence $a(\zeta)=\eta$
by the choice of $\eta$.
Since the second projection $p_2:G\times Y\rightarrow Y$ is flat, 
$p_2(\zeta)=y$.
Set $A:=\Cal O_{G\times X,\zeta}$.
Let $\varphi:\Cal O_{X,y}\rightarrow A$ and
$\psi:\Cal O_{X,\eta}\rightarrow A$ be the homomorphisms induced by
$p_2:G\times X\rightarrow X$ and $a:G\times X\rightarrow X$, respectively.
As $p_2$ and $a$ are finite-type flat local complete intersection morphisms 
\cite[(31.14)]{ETI},
$\varphi$ and $\psi$ are essentially of finite type flat homomorphisms
with complete intersection fibers.
As $\zeta$ is the generic point of 
a component of 
$p_2^{-1}(y)$, it is easy to see that
$\dim A=\dim \Cal O_{X,y}$.
Moreover, 
\[
A\otimes_{\O_{X,\eta}}\L_\eta\cong (a^*\L)_\zeta
\cong (p_2^*\L)_\zeta \cong A\otimes_{\O_{X,y}}\L_y.
\]
The last assertion is trivial.
\qed

\begin{corollary}
Let $S=\Spec k$, with $k$ a perfect field, and let $G$ be of finite type
over $S$.
Let $y$, $Y$, and $\eta$ be as in the theorem.
Then there exist $A$, $\varphi$ and $\psi$ as in the theorem
such that $\varphi$ and $\psi$ are regular.
\end{corollary}

\proof Replacing $G$ by $G\red$ if necessary, we may assume that $G$ is
$k$-smooth.
The assertion follows immediately by the theorem.
\qed

\begin{corollary}\label{eta-y.thm}
Let $y$ and $\eta$ be as in the theorem.
Then
$\dim \O_{X,y}\geq \dim \O_{X,\eta}$.
\end{corollary}

\proof $\dim \O_{X,y}=\dim A\geq \dim \O_{X,\eta}$.
\qed

\begin{corollary}\label{independence.thm}
Let $Y$ be as in the theorem, and $\eta_1$ and $\eta_2$ be the generic points
of irreducible components of $Y^*$.
Then $\dim\O_{X,\eta_1}=\dim\O_{X,\eta_2}$.
There are a noetherian local ring $A$ such that $\dim A=\dim \O_{X,\eta_1}$ 
and flat l.c.i.\ local homomorphisms $\varphi_i: \O_{X,\eta_i}\rightarrow A$
essentially of finite type such that for any quasi-coherent 
$(G,\O_X)$-module $\L$, $A\otimes_{\O_{X,\eta_1}}\L_{\eta_1}
\cong A\otimes_{\O_{X,\eta_2}}\L_{\eta_2}$.
If, moreover,  
$p_2:G\times X\rightarrow X$ 
is smooth, or $S=\Spec k$ with $k$ a perfect field and $G$ is of
finite type over $S$, then $\varphi_i$ can be taken to be regular.
\end{corollary}

\proof Let $\fp$ be the defining ideal of $Y$.
Let $\fq_i$ be the defining ideal of $Z_i$, where $Z_i$ is the closed
integral subscheme of $X$ whose generic point is $\eta_i$, for $i=1,2$.
By Corollary~\ref{Min-Min.thm}, $\fq_i^*\in\Min_G(\O_X/\fp^*)=\{\fp^*\}$.
Applying Corollary~\ref{eta-y.thm} to $y=\eta_1$ and $\eta=\eta_2$, 
$\dim \O_{X,\eta_1}\geq \dim \O_{X,\eta_2}$.
Similarly, $\dim \O_{X,\eta_2}\geq \dim \O_{X,\eta_1}$, and hence
$\dim \O_{X,\eta_1}= \dim \O_{X,\eta_2}$.
The other assertions are clear by Theorem~\ref{M-R-orig.thm}.
\qed

\begin{corollary}\label{M-R-general.thm}
Let $\Cal C$ and $\Cal D$ be classes of noetherian local rings, 
and $\Bbb P(A,M)$ a property of a finite module $M$ over a noetherian
local ring $A$.
Assume that:
\begin{description}
\item[(i)] If $A\in\Cal C$, $M$ a finite $A$-module with $\Bbb P(A,M)$, 
and $A\rightarrow B$ a flat l.c.i.\ \(resp.\ regular\) 
local homomorphism essentially of finite type,
then $B\in \Cal D$ and $\Bbb P(B,B\otimes_A M)$ holds.
\item[(ii)] If $A\rightarrow B$ is a flat l.c.i.\ \(resp.\ regular\) 
local homomorphism essentially of
finite type of noetherian local rings,
$M$ is a finite $A$-module,
and if $B\in \Cal D$ and $\Bbb P(B,B\otimes_AM)$ holds, 
then $A\in\Cal D$ and
$\Bbb P(A,M)$ holds.
\end{description}
If $\O_{X,\eta}\in\Cal C$ and $\Bbb P(\O_{X,\eta},\M_\eta)$ holds
\(resp.\ $\Bbb P(\O_{X,\eta},\M_\eta)$ holds and either 
$p_2:G\times X\rightarrow X$ is 
smooth, or
$S=\Spec k$ with $k$ a perfect field and $G$ is of finite type over $S$\),
then $\O_{X,y}\in\Cal D$ and $\Bbb P(\O_{X,y},\M_y)$ holds.
Conversely, if $\O_{X,y}\in\Cal C$ and $\Bbb P(\O_{X,y},\M_y)$ holds
\(resp.\ $\Bbb P(\O_{X,y},\M_y)$ holds and either $p_2:G\times X\rightarrow
X$ is smooth, or
$S=\Spec k$ with $k$ a perfect field and $G$ is of finite type over $S$\),
then $\O_{X,\eta}\in\Cal D$ and $\Bbb P(\O_{X,\eta},\M_\eta)$ holds.
\end{corollary}

\proof Because of $\psi:\O_{X,\eta}\rightarrow A$, 
$A\in\Cal D$ and $\Bbb P(A,A\otimes_{\O_{X,\eta}}\M_\eta)$ holds, by the
condition {\bf (i)}.
Since $A\otimes_{\O_{X,\eta}}\M_\eta\cong A\otimes_{\O_{X,y}}\M_y$
and we have $\varphi:\O_{X,y}\rightarrow A$, 
$\O_{X,y}\in\Cal D$ and $\Bbb P(\O_{X,y},\M_y)$ by the condition {\bf (ii)}.

The proof of the converse is similar.
\qed

\begin{corollary}\label{M-R.thm}
Let $y$ and $\eta$ be as in the theorem, and
$m$, $n$, and $g$ be nonnegative integers or $\infty$.
Then
\begin{description}
\item[(i)] If $\M_\eta$ is maximal Cohen--Macaulay
\(resp.\ of finite injecctive dimension, projective dimension $m$,
$\dim-\depth=n$, torsionless, reflexive, G-dimension $g$, zero\) as an
$\O_{X,\eta}$-module if and only if $\M_y$ is so as an
$\O_{X,y}$-module.
\item[(ii)] If $\O_{X,\eta}$ is a complete intersection, 
then so is $\O_{X,y}$.
\item[(iii)] Assume that $p_2:G\times X\rightarrow X$ 
is smooth,
or $S=\Spec k$ with $k$ a perfect field and $G$ is of finite type over $S$.
If $\O_{X,\eta}$ is regular, then so is $\O_{X,y}$.
\item[(iv)] Assume that $p_2:G\times X\rightarrow X$ is smooth,
or $S=\Spec k$ with $k$ a perfect field and $G$ is of finite type over $S$.
Let $p$ be a prime number, and assume that $\O_{X,\eta}$ is excellent.
If \mt{\O_{X,\eta}} is weakly $F$-regular \(resp.\ $F$-regular, $F$-rational\)
of characteristic $p$,
then so is $\O_{X,y}$.
If $\O_{X,y}$ is excellent and weakly $F$-regular of characteristic $p$, 
then $\O_{X,\eta}$ is
weakly $F$-regular of characteristic $p$.
\end{description}
\end{corollary}

\proof {\bf (i)} 
Let $\Cal C=\Cal D$ be the class of all noetherian local rings,
and $\Bbb 
P(A,M)$ be the property ``$M$ is a maximal Cohen--Macaulay $A$-module''
in Corollary~\ref{M-R-general.thm}.
The assertion follows immediately by Corollary~\ref{M-R-general.thm}.
Similarly for other properties.

{\bf (ii)} 
Let $\Cal C=\Cal D$ be the class of complete intersection noetherian
local rings, and $\Bbb P(A,M)$ be ``always true.''

{\bf (iii)} Let $\Cal C=\Cal D$ be the class of regular local rings.

{\bf (iv)} Let $\Cal C$ be the class of excellent weakly $F$-regular
local rings of characteristic $p$,
and $\Cal D$ be the class of weakly $F$-regular local rings of
characteristic $p$.
By \cite{HH}, {\bf (i)} and {\bf (ii)} 
of Corollary~\ref{M-R-general.thm} holds.
Similarly for $F$-regularity, see \cite{HH}.
For $F$-rationality, see \cite{Velez}.
\qed

\begin{corollary}\label{independence2.thm}
Let $Y$ be a $G$-primary closed subscheme of $X$.
Let $\eta_1$ and $\eta_2$ be the generic points of irreducible components of
$Y$.
Let $d$, $\delta$, $m$, and $g$ be nonnegative 
integers 
or $\infty$.
\begin{description}
\item[(i)] $\M_{\eta_1}$ is maximal Cohen--Macaulay
\(resp.\ of finite injective dimension, dimension $d$, depth $\delta$, 
projective dimension $m$, torsionless, reflexive, G-dimension $g$, zero\) as
an $\O_{X,\eta_1}$-module if and only if
the same is true of $\M_{\eta_2}$ as an $\O_{X,\eta_2}$-module.
\item[(ii)] $\O_{X,\eta_1}$ is a complete intersection if and only if
$\O_{X,\eta_2}$ is so.
\item[(iii)] 
Assume that $p_2:G\times X\rightarrow X$ 
is smooth, or $S=\Spec k$ with $k$ a perfect field
and $G$ is of finite type over $S$.
Then $\O_{X,\eta_1}$ is regular if and only if $\O_{X,\eta_2}$ is so.
Assume further that $X$ is locally excellent \(that is, the 
all local rings of $X$ are excellent\).
Then $\O_{X,\eta_1}$ is of characteristic $p$ and 
weakly $F$-regular \(resp.\ $F$-regular, $F$-rational\) 
if and only if $\O_{X,\eta_2}$ is so.
\end{description}
\end{corollary}

\proof Follows immediately from Corollary~\ref{independence.thm}.
\qed

\paragraph
Let $Y=V(\Q)$ be as in the corollary.
Then we say that $\M$ is maximal Cohen--Macaulay (resp.\ of finite injective
dimension, dimension $d$, depth $\delta$, projective dimension $m$, 
torsionless, reflexive, $G$-dimension $g$) along $Y$ (or at $\Q$) if
$\M_\eta$ is so for the generic point $\eta$ of some (or equivalently, any) 
irreducible component of $Y$.
We say that $X$ is complete intersection along $Y$ (or at $\Q$) 
if $\O_{X,\eta}$ is
a complete intersection for some (or equivalently, any) $\eta$.
Assume that $p_2:G\times X\rightarrow X$ 
is smooth, or $S=\Spec k$ with $k$ a perfect field
and $G$ is of finite type over $S$.
We say that $X$ is regular along $Y$ (or at $\Q$) 
if $\O_{X,\eta}$ is a regular local
ring for some (or equivalently, any) $\eta$.
Assume further that $X$ is a locally excellent $\Bbb F_p$-scheme.
We say that $X$ is weakly $F$-regular (resp.\ $F$-regular, $F$-rational) 
along $Y$ (or at $\Q$) 
if $\O_{X,\eta}$ is so for some (or equivalently, any) $\eta$.

\begin{corollary}\label{M-R-graded.thm}
Let $p$ be a prime number, and
$A$ a $\Bbb Z^n$-graded noetherian ring.
Let $P$ be a prime ideal of $A$, and $P^*$ be the prime ideal of $A$
generated by the homogeneous elements of $P$.
If $A_{P^*}$ is excellent of characteristic $p$ and is
weakly $F$-regular \(resp.\ $F$-regular, $F$-rational\), then
$A_P$ is weakly $F$-regular \(resp.\ $F$-regular, $F$-rational\).
If $A_P$ is excellent of characteristic $p$ and is weakly $F$-regular,
then $A_{P^*}$ is weakly $F$-regular.
\end{corollary}

\proof Let $S=\Spec \Bbb Z$, $G=\Bbb G_m^n$, and $X=\Spec A$.
If $y=P$, then $\eta$ in the theorem is $P^*$.
The assertion follows immediately by Corollary~\ref{M-R.thm}.
\qed

\def\citinfo{Cf.\ \cite[(4.7)]{HH2}}
\begin{corollary}[\citinfo]\label{graded-enough.thm}
Let $A$ be a $\Bbb Z^n$-graded locally excellent ring of characteristic $p$.
If $A_\fm$ is $F$-regular \(resp.\ $F$-rational\) for any graded maximal
ideals \(that is, $G$-maximal ideal for $G=\Bbb G_m^n$\), 
then $A$ is $F$-regular \(resp.\ $F$-rational\).
\end{corollary}

\proof If $Q$ is a graded prime ideal, then it is contained in a graded
maximal ideal.
So $A_Q$ is $F$-regular (resp.\ $F$-rational, see \cite[(4.2)]{HH}).

Now consider any prime ideal $\fp$ of $A$.
Then $A_{\fp^*}$ is $F$-regular (resp.\ $F$-rational), since 
$\fp^*$ is a graded prime ideal.
Hence $A_{\fp}$ is so by Corollary~\ref{M-R-graded.thm}.
So $A$ is $F$-regular (resp.\  $F$-rational).
\qed

\begin{remark}
Let $k$ be an ($F$-finite) field of characteristic $p$, 
and $A=\bigoplus_{n\geq 0}A_n$, $A_0=k$,
a positively graded finitely generated $k$-algebra.
Set $\fm:=\bigoplus_{n>0}A_n$.
Lyubeznik and Smith \cite{LS} proved that 
if $A_\fm$ is weakly $F$-regular, then $A$ is (strongly) $F$-regular.
\end{remark}

\begin{corollary}\label{graded-deformation.thm}
Let $A=\bigoplus_{n\geq 0}A_n$ be an $\Bbb N$-graded 
locally excellent noetherian
ring of characteristic $p$.
Let $t\in A_+:=\bigoplus_{n>0}A_n$ be a nonzerodivisor of $A$.
If $A/tA$ is $F$-rational, then $A$ is $F$-rational.
\end{corollary}

\proof Note that $A/tA$ is Cohen--Macaulay \cite[Proposition~0.10]{Velez}.
By Corollary~\ref{graded-enough.thm}, 
it suffices to show that $A_{\fm}$ is $F$-rational for any
graded maximal ideal (i.e., $G$-maximal ideal for $G=\Bbb G_m^n$).
It is easy to see that $\fm=\fm_0+A_+$.
Hence $t\in \fm$.
Since $A_\fm/tA_\fm$ is Cohen--Macaulay $F$-rational, $A_\fm$ is also
$F$-rational by \cite[(4.2)]{HH}.
\qed

\begin{corollary}\label{filtration-deformation.thm}
Let $A$ be a ring of characteristic $p$, and $(F_n)_{n\geq 0}$ 
a filtration of $A$.
That is, each $F_n$ is an additive subgroup of $A$, 
$F_0\subset F_1\subset F_2\subset\cdots\subset A$, 
$1\in F_0$, $F_iF_j\subset F_{i+j}$, and $\bigcup_{n\geq 0}F_n=A$.
Set $R=\bigoplus_{n\geq 0}F_n t^n\subset A[t]$, and $G=R/tR$.
If $G$ is \(locally\) excellent noetherian and $F$-rational, then
$A$ is also \(locally\) excellent noetherian and $F$-rational.
\end{corollary}

\proof There exist homogeneous elements $a_1 t^{n_1},\cdots ,a_r t^{n_r}$ of
$R$ such that their images in $G$ generate 
the ideal $G_+$, the irrelevant ideal of $G$.
Then $R_+$ is generated by $t, a_1t^{n_1}, \ldots,a_rt^{n_r}$.
Since $F_0\cong G/G_+$, $F_0$ is (locally) excellent noetherian.
Since $R$ is generated by $t, a_1t^{n_1}, \ldots,a_rt^{n_r}$ as an 
$F_0$-algebra, $R$ is also (locally) excellent noetherian.
Then by Corollary~\ref{graded-deformation.thm}, 
$R$ is $F$-rational.
Hence $R[t^{-1}]\cong A[t,t^{-1}]$ is (locally) excellent and $F$-rational.
Hence $A$ is (locally) excellent and $F$-rational.
\qed

\end{document}